\documentclass[reqno,10pt]{article}
\usepackage{amsmath, latexsym, amsfonts, amssymb, amsthm}
\usepackage{graphicx, color,hyperref,dsfont,tikz,caption,wrapfig,subfig}
\usepackage{makeidx}
\usepackage{pstricks}

\hypersetup{
    linktoc=page,
    linkcolor=red,          
    citecolor=blue,        
    filecolor=blue,      
    urlcolor=cyan,
    colorlinks=true           
}

\numberwithin{equation}{section}

\newtheorem{theorem}{Theorem}[section]
\newtheorem{lemma}[theorem]{Lemma}

\newtheorem{rem}[theorem]{Remark}

\newtheorem{conjecture}[theorem]{Conjecture}

\makeindex

\renewcommand{\tilde}{\widetilde}          
\DeclareMathSymbol{\leqslant}{\mathalpha}{AMSa}{"36} 
\DeclareMathSymbol{\geqslant}{\mathalpha}{AMSa}{"3E} 
\DeclareMathSymbol{\eset}{\mathalpha}{AMSb}{"3F}     
\renewcommand{\leq}{\;\leqslant\;}                   
\renewcommand{\geq}{\;\geqslant\;}                   
\renewcommand{\epsilon}{\varepsilon}

\newcommand{\C}{\mathbb{C}}

\newcommand{\R}{\mathbb{R}}
\newcommand{\Z}{\mathbb{Z}}
\newcommand{\N}{\mathbb{N}}

\newcommand{\E}{\mathds{E}}

\newcommand{\gb}{\beta}
\newcommand{\gep}{\varepsilon}

\def\P{{\mathbb P}}

\def\N{{\mathbb N}}

\def\1{{\mathds{1}}}

\def\ee{\mathrm{e}}

\title{The glassy phase of complex branching Brownian motion}

\begin{document}

\maketitle
\begin{center}

{ Thomas Madaule \footnotemark[1],\ R\'emi Rhodes \footnotemark[2],\ 
Vincent Vargas \footnotemark[3]}

\bigskip

 \footnotetext[1]{Universit{\'e} Paris-13} 
\footnotetext[2]{Universit{\'e} Paris-Dauphine, Ceremade, F-75016 Paris, France. Partially supported by grant ANR-11-JCJC  CHAMU}
\footnotetext[3]{Ecole normale sup\'erieure, DMA. Partially supported by grant ANR-11-JCJC  CHAMU}

\end{center}

\begin{abstract}
In this paper, we study complex valued branching Brownian motion in the so-called glassy phase, or also called phase II. In this context, we prove  a limit theorem for the complex partition function hence confirming a conjecture formulated by Lacoin and the last two authors in a previous paper on complex Gaussian multiplicative chaos. We will show that the limiting partition function can be expressed as a product of a Gaussian random variable, mainly due to the windings of the phase, and a stable transform of   the so called derivative martingale, mainly due to the clustering of the modulus. The proof relies on the fine description  of the extremal process available in the branching Brownian motion context. 
\end{abstract}
\vspace{0.3cm}
\footnotesize


\noindent{\bf Key words or phrases:} branching Brownian motion, freezing, glassy phase.


\normalsize

\tableofcontents

\section{Introduction}

\begin{figure}[t]
\centering
\begin{tikzpicture}[scale=0.8]
\shade[left color=green!10,right color=green] (3,0) -- (12,0) -- (12,6) -- (3,6) ;
\shade[top color=red, bottom color=red!20] (3.5,3.5) -- (12,3.5) -- (12,10) -- (0,10) -- (0,7);
\shade[left color=yellow,right color=yellow!10] (5,0) arc (0:45:5) -- (0,7) -- (0,0) -- (5,0);
\draw [line width=2.5pt,color=blue](3.53,3.5) -- (12,3.5);
\draw [line width=2.5pt,color=blue](3.6,3.5) -- (0,7);
\draw [line width=2.5pt,color=blue](5,0)  arc (0:45:5);
\draw[style=dashed,line width=2pt,->] (0,0) -- (0,8.2) node[left]{{\Large $\gamma$}};
\draw[style=dashed,line width=2pt,->] (0,0) -- (12.2,0)node[below]{{\Large $\beta$}};
\draw (7,3.5)  node[above]{{\large $\gamma=\sqrt{\frac{d}{2 }}$}};
\draw (2,5.2)  node[above,rotate=-45,line width=2pt]{{\large $\gamma+\beta=\sqrt{2d}$}};
\draw (4.8,1.5)  node[above,rotate=-70,line width=2pt]{{\large $\gamma^2+\beta^2=d$}};
\draw (1,2)  node[right]{{\Large Phase I}};
\draw (5,7)  node[right]{{\Large Phase II}};
\draw (7,2)  node[right]{{\Large Phase III}};
\end{tikzpicture}
\caption{Phase diagram}
\label{diagram}
\end{figure}
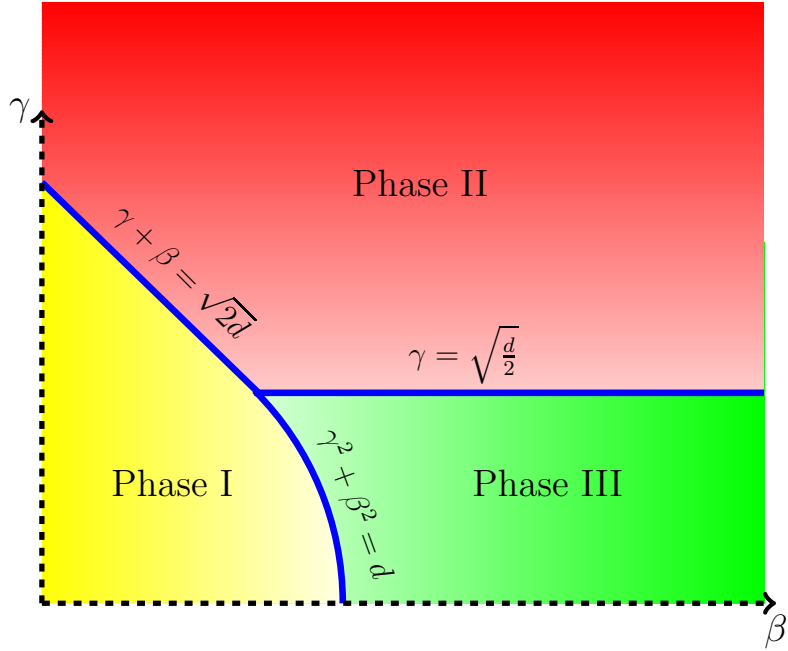

In a recent article \cite{LRV}, the authors studied complex Gaussian multiplicative chaos, a complex extension of classical Gaussian multiplicative chaos (see \cite{review} for a review on Gaussian multiplicative chaos). More precisely, consider two independent and logarithmically correlated Gaussian fields $X,Y$ on a subdomain $\Omega \subset \R^d$
\begin{equation*}
\E[X(x)X(y)] = \E[Y(x)Y(y)] \underset{|y-x| \to 0}{\sim}  \ln \frac{1}{|y-x|}.
\end{equation*}
We denote $\mathcal{D}(\Omega)$ the space of smooth functions with compact support in $\Omega$ and $\mathcal{D}(\Omega)$ the space of distributions (in the sense of Schwartz). They adressed the problem of finding a proper renormalization as well as the limit of the family of complex random distributions
\begin{equation}\label{mesuresapprox}
M^{\gamma,\beta}_\varepsilon( \varphi)=\int_{\Omega} e^{\gamma X_\epsilon(x)+i\beta Y_\epsilon(x)} \varphi(x)\,dx,  \quad \varphi \in \mathcal{D}(\Omega)
\end{equation}
where $X_\varepsilon, Y_\varepsilon$ are appropriate regularizations (say of variance of order $\ln \frac{1}{\varepsilon}$) which converge to $X,Y$ and $\gamma,\beta$ are real constants. In this setting, they recovered the phase diagram of figure \ref{diagram} which was first discovered in the pioneering work \cite{derrida} in the simpler context of discrete multiplicative cascades, i.e. when $X_\varepsilon,Y_\varepsilon$ are independent branching random walks on a tree-like structure. More precisely, the authors of \cite{derrida} computed the free energy of the total mass (or partition function) 
\begin{equation*}
\underset{\varepsilon \to 0} {\lim}  \: \frac{1}{\ln \varepsilon } \ln |  M^{\gamma,\beta}_\varepsilon(\Omega) |
\end{equation*}
for a subset $\Omega$ and found phase transitions according to a diagram similar to our figure \ref{diagram}. In particular, they distinguished three phases I, II and III which we have indicated on the figure. The work \cite{LRV} is a step further in understanding the limit of \eqref{mesuresapprox}. Indeed, the framework of \cite{LRV} is that of finding a deterministic sequence $c(\varepsilon)$ such that $c(\varepsilon) M^{\gamma,\beta}_\varepsilon $ converges to a non trivial limit in the space of distributions (see also the interesting and related works \cite{bar:comp1,bar:comp2}). In a series of works \cite{Rnew7,Rnew12,LRV,MRV}, this question was essentially solved for the phases I and III (and their frontiers) but left unanswered in phase II. However, it was conjectured that in phase II, the behaviour of $M^{\gamma,\beta}_\varepsilon$ is mainly ruled by two phenomena: the local intensity of this complex measure is dominated by the local maxima of the field $X_\varepsilon$ whereas the overall phase resulting from  the (strong) windings of the field $Y_\epsilon$ asymptotically behaves like a white noise. This led to the following freezing conjecture corresponding to the so-called glassy phase (the freezing and glassy phase terminology comes from physics, see \cite{CarDou,Fyo,rosso} for example):

 \begin{conjecture}\label{conjII}
 Let $\beta>0$ and $\gamma>\sqrt{\frac{d}{2}}$ be such that $\beta>\max(\sqrt{2d}-\gamma,0)$. Set $\alpha=\sqrt{\frac{d}{2}}\frac{1}{\gamma}$. There exists some constant $\sigma:=\sigma(\gamma,\beta)>0$ such that we get the following convergence in law:
\begin{equation}\label{eqconjecture}
\left( (\ln \frac{1}{\gep})^{\frac{3 \gamma}{2\sqrt{2d}}} \gep^{   \gamma \sqrt{2d}-d} M^{\gamma,\gb}_\gep(A)\right)_{A\subset \R^d} \Rightarrow \left(\sigma W_{ N^\alpha_{M'}}(A)\right)_{A\subset \R^d},   
 \end{equation}
 where, conditionally on  $N^\alpha_{M'}$,  $W_{ N^\alpha_{M'}}$ is a complex Gaussian random measure with intensity  $ N^\alpha_{M'}$ and  $N^\alpha_{M'}$ is a $\alpha$-stable random measure with intensity $M'$, namely a  random distribution  whose law is characterized by $\E[e^{-q W_{N^\alpha_{M'}(A)}} ]=\E[e^{-q^{2\alpha} M'(A)}]$ for every $q\geq 0$ and every bounded Borelian set $A$ .  
 \end{conjecture}

Let us finally mention that  a  result similar to \eqref{eqconjecture} is proved  in the paper \cite{MRV} in the real case, i.e. on the frontier of phase II hence for $\gamma>\sqrt{2d}$ and $\beta=0$ (see also \cite{arguin, arguin2, Louisdor} for related results). Analogous results in the real case were also derived recently for the Branching Random Walk in \cite{BKNSW,Rnew3,madaule,Webb}. Recall that in the context of the real Branching Random Walk, these problems have received much attention since the works \cite{DS,mandelbrotstar}.

The purpose of this work is to prove the analogue of conjecture \eqref{conjII} in the context of the simpler but related model, the so-called branching Brownian motion (BBM) where the approximations $X_\varepsilon,Y_\varepsilon$ are defined by particles which split along a Poisson process and then perform independent Brownian motions: see the next section for precise definitions. Let us mention that, up to some technical adaptations, it should be possible to prove in the BBM context results analogue to \cite{LRV} and in particular to recover a phase diagram similar to figure \ref{diagram}. Over the past years, there has been impressive progress on the the study of BBM since the seminal works \cite{B1,B2,LS}: this progess has culminated in the works \cite{ABBS,ABK1,ABK2,ABK3}. Thanks to these achievements, it is possible to know with high precision the behaviour of the extreme particles of the BBM which dominate phase II. Though our work in the context of BBM relies on the fine results of \cite{ABBS,ABK1,ABK2,ABK3}, we believe that it gives insights on the mechanism involved behind the conjectured convergence \eqref{eqconjecture}: this will be discussed in more detail in the next section.         

The paper is organized as follows. In the next section, we define the setup and cite the main result of the paper namely theorem \ref{maintheorem}. We also include a discussion on related models, like the branching random walk or the logarithmically correlated Gaussian fields considered in \cite{MRV}. Special emphasis will be given to the case of the maximum of the discrete Gaussian Free Field which has received a lot of attention recently \cite{arguin2,Louisdor,DZ,Rnew7,Rnew12}. In the following section, we prove theorem \ref{maintheorem}.

\section{Setup and main result}

\subsection{Setup and main result}
In this paper, we will study the branching Brownian motion (BBM for short). Start with a single particle which performs standard Brownian motion starting from $0$ up to an exponential holding time $T$ with parameter $\lambda=1$. At time $T$, the particle splits into two new particles, both of them following a new independent copy of the same process starting from its place of birth. Both new particles thus move according to a standard Brownian motion and split into two after an exponential holding time and so on.
We introduce $N(t)$ the associated Poisson  point process which counts the number of particles at time $t$ and $(\bar{X}_{i}(t))_{1 \leq i \leq N(t)}$ the (increasingly ordered) positions of the particles.

We then introduce the properly normalized and shifted  quantity
\begin{equation*} 
X_{i}(t)= \sqrt{2} \bar{X}_i(t)+2t,
\end{equation*}
in order to have:
\begin{equation}
\label{critical} \E \left[ \sum_{i=1}^{N(t)} e^{-X_i(t)} \right] =1,\qquad \E\left[\sum_{i=1}^{N(t)} X_i(t) e^{-X_i(t)} \right] = 0,\qquad \forall t>0.
\end{equation}

On the same probability space, we consider particles which split according to the \textbf{same} Poisson point process $N(t)$ but follow Brownian motions that are independent of those involved in the definition of $X$.  We consider $(\bar{Y}_{i}(t))_{1 \leq i \leq N(t)}$ the positions of these new particles.

We introduce the random measure
\begin{equation} \label{defmesure}
\mathcal{N}_{t}(dX,dY)= \sum_{i=1}^{N(t)} \delta_{(X_i(t),Y_i(t))}
\end{equation}
We will also consider the measure  $\bar{\mathcal{N}}_{x,t}$ which corresponds to   the measure $\mathcal{N}_{t}$ conditioned to the event that all particles $(X_i(t))_{1 \leq i \leq N(t)}$ are above $-x$. If $f$ is some continuous function, we denote
\begin{equation*}
<f(X,Y), \mathcal{N}_{t}(dX,dY)> := \sum_{i=1}^{N(t)} f(X_i(t),Y_i(t))
\end{equation*}     
and similarly for $\bar{\mathcal{N}}_{x,t}$.

In order to state our results, we introduce the limit of the derivative martingale $M'$ given by the following limit (first derived in \cite{neveu})
\begin{equation*}
M':= \underset{t \to \infty}{\lim} \sum_{i=1}^{N(t)}  X_i(t) e^{-X_i(t)}
\end{equation*}
Recall the following classical convergence in law of the minimum obtained in \cite{B2}
\begin{equation}\label{convBram}
X_1(t)-\frac{3}{2} \ln t \underset{t \to \infty}{\rightarrow}  W
\end{equation}
where $W$ is some random variable satisfying $\P(W \leq x ) \underset{x \to - \infty}{\sim} c_\star |x|e^{x}$ and $c_\star$ some constant.

We are interested in studying the variable
\begin{equation}\label{defvariable}
 t^{\frac{3 \gamma}{2}}  \sum_{i=1}^{N(t)}e^{-\gamma X_{i}(t)+\beta \sqrt{2}\bar{Y}_{i}(t) }
\end{equation}
in the so-called phase II, i.e. $\beta>(1-\gamma)_+$ and $\gamma>\frac{1}{2}$ where for a real $x$ we set $x_+= \max(x,0)$. To state our main result, we recall that a random variable $\mathcal{G}$ is a standard complex Gaussian random variable if $\mathcal{G}= \mathcal{G}_1+i \mathcal{G}_2$ where $\mathcal{G}_1,\mathcal{G}_2$ are two independent standard real Gaussian variables. The following theorem is the main result of the paper:

\begin{theorem}\label{maintheorem}
For $(\gamma,\beta)$ in phase II, there is some constant $c(\gamma,\beta)>0$ such that we have the following convergence in law
\begin{equation}\label{eqtheo}
 t^{\frac{3 \gamma}{2}} \sum_{i=1}^{N(t)}e^{-\gamma X_{i}(t)+i \sqrt{2} \beta \bar{Y}_{i}(t) }  \underset{t \to \infty} {\rightarrow} c(\gamma,\beta)   N_\alpha^{1/2} \mathcal{G}  
 \end{equation}
 where $\mathcal{G}$ is a standard complex Gaussian random variable independent from $N_\alpha$ which is a $\alpha$-stable random variable with intensity $M'$ and $\alpha=\frac{1}{2 \gamma}$. More precisely, the law of $N_\alpha$ is characterized by 
\begin{equation*} 
\E[e^{-q N_\alpha }  ]=\E[e^{-q^\alpha M'}].
\end{equation*} 
for all $q\geq 0$. 
\end{theorem}

\begin{rem}
In a recent work \cite{KaKli}, the authors showed a result similar to \eqref{eqtheo} in the context of the complex REM, i.e. when the variables $(X_i(t),Y_i(t))_i$ form an iid sequence of centered Gaussian vectors. In the REM context, one must replace the renormalization $ t^{\frac{3 \gamma}{2}}$ in \eqref{eqtheo} by $t^{\frac{\gamma}{2}}$ and the limiting law is of the same form as the right hand side of \eqref{eqtheo} with the variable $N_\alpha$ distributed as a standard stable distribution (whereas, in the BBM context of theorem \ref{maintheorem}, the variable $N_\alpha$ is stable conditionally to $M'$). Note that in our case the result is different because of the strong correlations in the BBM; in particular, the methods of \cite{KaKli} which rely on the summation theory of triangular arrays of independent random vectors can not be adapted here.  
\end{rem}

\subsection{Heuristic of the proof}
In this subsection, we start by giving an insight on the proof of theorem \ref{maintheorem} which will enable us to discuss other related models: the branching random walk and the discrete GFF.
First, introduce the set $N_{loc}(t)$ of local minima of $X_{i}(t)$ that are close to $\frac{3}{2} \ln t$, i.e. those particles which are at distance of order 1 from $\frac{3}{2} \ln t$ and that are smaller than    all the particles sharing with them a common ancestor at distance of order 1. In phase II, the variable \ref{defvariable} concentrates on the local minima along with the close neighbours that do not have atypical high values, which constitute the so-called decoration. Therefore, the variable \eqref{defvariable} is roughly equal for large $t$ to 
\begin{equation*}
  t^{\frac{3 \gamma}{2}} \sum_{u \in N_{loc}(t)} \sum_{ j \approx u, X_{j}(t) \approx X_{u}(t)} e^{-\gamma X_{j}(t)+i \sqrt{2}\beta\bar{Y}_{j}(t) }
\end{equation*}
where $x \approx y$ means that $|x-y|$ is of order 1. Now, one can rewrite the above quantity in the following way
\begin{equation*}
 t^{\frac{3 \gamma}{2}} \sum_{u \in N_{loc}(t)} e^{-\gamma X_{u}(t) }  \:\left ( e^{i \sqrt{2}\beta \bar{Y}_{u}(t) }  \sum_{ j \approx u, X_{j}(t) \approx X_{u}(t)} e^{-\gamma (X_{j}(t)-X_{u}(t))+i \sqrt{2}\beta (\bar{Y}_{j}(t)-\bar{Y}_{u}(t) ) }  \right )
\end{equation*}

From the results of \cite{ABBS,ABK1,ABK2,ABK3}, the sum $ t^{\frac{3 \gamma}{2}} \sum_{u \in N_{loc}(t)} e^{-\gamma X_{u}(t) }$ converges to $\sum_{u \geq 1} e^{- \gamma \Delta_u } $ where $(\Delta_u)_{u \geq 1}$ is a Poisson point process (PPP) with intensity $cM' \,e^x dx$ where $c>0$ is some constant. 

Since the local minima are far apart, each sum $( e^{i \sqrt{2} \beta \bar{Y}^{(u)}(t)   }  \sum_{ j \approx u, X_{j}(t) \approx X_{u}(t)}  \cdots)_u$ is asymptotically independent for different values of $u$. From the results of \cite{ABBS,ABK1,ABK2,ABK3}, one can also deduce that each term 
\begin{equation*}
\sum_{ j \approx u, X_{j}(t) \approx X_{u}(t)} e^{-\gamma (X_{j}(t)-X_{u}(t))+i \sqrt{2}\beta(\bar{Y}_{j}(t)-\bar{Y}_{u}(t)} )
\end{equation*}
converges in law to some non trivial variable $Z^{(u)}$ (which is painful to describe).

Finally, if $\mathcal{N}$ is a standard Gaussian the variable $e^{i\alpha \mathcal{N}}$ converges in law as $\alpha \to \infty$ to a random variable uniformly distributed on the unit circle. 
Hence $e^{ i \sqrt{2} \beta \bar{Y}_{u}(t) }$ converges in law to a variable $U_u$ uniformly distributed on the unit circle and independent from $Y^{(u)}$. Gathering the above considerations, we see that the variable \eqref{defvariable} converges to 
\begin{equation*}
\sum_{u \geq 1} e^{- \gamma \Delta_u } U_u Z^{(u)}
\end{equation*}
where $(U_uZ^{(u)})_{u \geq 1}$ is an i.i.d. sequence of \textbf{isotropic} random variables. Though we do not have a friendly description of the variable $U_uZ^{(u)}$, the scaling property of the Poisson sequence and the isotropy of $U_uZ^{(u)}$ yield representation \eqref{eqtheo}. In fact, a similar mechanism is behind the freezing phenomenon in the real case $(\gamma>1, \beta=0)$; indeed, in this case, the i.i.d property of the decoration combined to the scaling property of the PPP yield a stable distribution.

\subsection{Discussion of other models}

In the case of the branching random walk, it should be possible to prove analogues of \cite{ABBS,ABK1,ABK2,ABK3} though it certainly requires non trivial technical difficulties to adapt the proofs of \cite{ABBS,ABK1,ABK2,ABK3}. Therefore, proving a result similar to theorem \ref{maintheorem} is clearly within reach (in the lognormal and even the non lattice case).

In the case of the discrete GFF, the situation is a bit more involved and does not just require technical adaptations: this is due to the fact that the correlations do not involve a hierarchical structure. Consider a discrete GFF $X_\varepsilon$ on a square grid of size $\varepsilon$ in a fixed bounded domain $D$ with the normalization $\E[ X_\varepsilon(x)^2]= 2 \ln \frac{1}{\varepsilon}+2 \ln C(x,D) +o(1) $ where $C(x,D)$ denotes the conformal radius of a point $x \in D$. Fix $l>0$. We introduce the set $N_{loc}^{(l)}(\varepsilon):= \lbrace x_1, ..., x_{J_\varepsilon^{(l)}}  \rbrace$ of coordinates of the local maxima of $X_\varepsilon$ that are in the interval 
\begin{equation*}
[2 \ln \frac{1}{\varepsilon}+ \frac{3}{2}\ln \ln \frac{1}{\varepsilon}-l,2 \ln \frac{1}{\varepsilon}+ \frac{3}{2}\ln \ln \frac{1}{\varepsilon}+l ]\end{equation*} 
In view of the results of \cite{Louisdor,MRV}, it is natural to conjecture that the following convergence in law holds for all $k$
\begin{align}
& \Big (   \big(  x_u , X(x_u),  \big(  X(y) - X(x_u) , \frac{y - x_u}{\varepsilon } \big)_{\lbrace y; \: |X(y) - X(x_u)| \leq k \: \text{and} \: |y - x_{\textcolor{red}{u}}|  \leq \varepsilon k   \rbrace}\big)  \Big )_{u \leq J_\varepsilon^{(l)}}  \nonumber   \\
& \underset{\varepsilon \to 0} {\rightarrow}  (z_u,\Delta_u, \mu_u  ( . 1_{|.|\leq k}),\rho_u ( . 1_{|.|\leq k})  )_{u \geq 1, \Delta_u \in [-l,l]} \label{convconj}
\end{align}  
 where $(z_u, \Delta_u)_{u \geq 1}$ is a Poisson Point process (PPP) with intensity $c \: C(x,D)^2 M'(dx)\times e^{-y} dy$ where $M'$ is the derivative martingale constructed in \cite{Rnew7,Rnew12}, $c$ some constant and $(\mu_u,\rho_u)_{u \geq 1}$ is an i.i.d. sequence of couples of point processes that are independent from the $(\Delta_u)_{u \geq 1}$. The law of $\rho_u$ should be isotropic. Recall the remarkable result of \cite{Louisdor} where the authors prove that $ (   (  x_u , X(x_u)  )_{u \leq J_\varepsilon^{(l)}}$ converges to a PPP with intensity $Z(dx)\times \1_{\{y \in [-l,l]\}}e^{-y} dy$ where $Z(dx)$ should coincide with $c \:C(x,D)^2 M'(dx)$.      
 
Nonetheless, in order to adapt our result to this context, one still has to reinforce the conjectured convergence \eqref{convconj} by adding information on the "decorrelation time" $\varepsilon$ of two points in the point process $\rho$. This is certainly a non trivial issue that requires a fine analysis of the discrete GFF.

\section{Proof of theorem \ref{maintheorem}}
We first recall the following useful lemma, the so-called many-to-one lemma, which states that for all nonnegative function $F$
\begin{align} 
\nonumber
\E[ \sum_{i=1}^{N(t)}  F((X_i(s))_{0 \leq s \leq t})   ] &=e^t \E[ F  ( ( \sqrt{2}B_s+ 2s)_{0 \leq s \leq t}) ] \\
& = \E[ e^{\sqrt{2}B_t} F  ( ( \sqrt{2}B_s)_{0 \leq s \leq t}) ] \label{manytoone}
\end{align}
where $(B_s)_{s \geq 0}$ is a standard Brownian motion.

\subsection{Proof of theorem \ref{maintheorem}}

Given the technical lemmas of the next subsection, it is not very difficult to conduct the proof. 
Let $F$ be some bounded and Lipschitz function from $\C$ to $\R$. We will additionally suppose that $F$ is bounded by $1$ and $1$-Lipschitz. 
By lemma \ref{complextension}, there exists $C_k$ such that $\underset{k \to \infty}{\lim}C_k=0$ and 
\begin{equation*}
\underset{t \to \infty}{\overline{\lim}}  \E \left[  \left |  F(t^{\frac{3 \gamma}{2}} \sum_{i=1}^{N(t)}e^{-\gamma X_{i}(t)+i \sqrt{2} \beta \bar{Y}_{i}(t) })- F(t^{\frac{3 \gamma}{2}} \sum_{i=1}^{N(t)} \1_{\{X_i(t) \leq \frac{3}{2}\ln t +k \}}e^{-\gamma X_{i}(t)+i \sqrt{2} \beta \bar{Y}_{i}(t) })  \right |  \right  ] \leq C_k
\end{equation*}

Following \cite{ABBS}, for $l \geq 0$, we introduce $\mathcal{H}_l(t)$ the set of particles which  are the first in their line of descent up to time $t$ to hit $l$. Since $X_1(t)$ converges almost surely to infinity as $t$ goes to infinity, $\mathcal{H}_l(t)$ is constant for $t$ (random) large enough and equal to a set that we will denote $\mathcal{H}_l$. Observe that $\mathcal{H}_l$ is finite almost surely. 

For each $u \in \mathcal{H}_l(t)$, we consider the ordered descendants $(X^{u}_i(t))_{1 \leq i \leq N^u(t)}$ up to time $t$. Then, we have 
\begin{align}   
 & t^{\frac{3 \gamma}{2}} \sum_{i=1}^{N(t)}e^{-\gamma X_{i}(t)+i \sqrt{2} \beta \bar{Y}_{i}(t) } \1_{\{X_i(t) \leq \frac{3}{2}\ln t +k \}} \nonumber \\
& = \sum_{u \in \mathcal{H}_l(t)}  t^{\frac{3 \gamma}{2}} e^{-\gamma X_1^u(t)+i \sqrt{2} \beta \bar{Y}_1^u(t)   }  \sum_{j=1}^{N^u(t)}\1_{\{X_j^u(t) \leq \frac{3}{2}\ln t +k \}}  e^{-\gamma (X_{j}^u(t)-X_{1}^u(t))+i \sqrt{2} \beta (\bar{Y}_{j}^u(t)-\bar{Y}_{1}^u(t)  ) } +A_{t,l,k} \label{truc}  
\end{align}
where $A_{t,l,k}$ corresponds to the sum on the $i$ which are not descendants of $u \in \mathcal{H}_l(t)$. Since $X_1(t)$ converges almost surely to infinity as $t$ goes to infinity, the variable $\sup_{k \geq 0} A_{t,l,k}$ converges 
almost surely to $0$ as $t$ goes to infinity. Hence, we just have to study the convergence of $\E[ F(\sum_{u \in \mathcal{H}_l(t)} \cdots)]$ where $\sum_{u \in \mathcal{H}_l(t)} \cdots$ is defined in equality \ref{truc} .
We introduce $\tau^u_{i,j}(t)$ the splitting times of particles $X_{i}^u(t)$ and $X_{j}^u(t)$. Now, we have
\begin{align}
\nonumber& \E \left[ F \left( \sum_{u \in \mathcal{H}_l(t)}  t^{\frac{3 \gamma}{2}} e^{-\gamma X_1^u(t)+i \sqrt{2} \beta \bar{Y}_1^u(t)   }  \sum_{j=1}^{N^u(t)}\1_{\{X_j^u(t) \leq \frac{3}{2}\ln t +k \}}  e^{-\gamma (X_{j}^u(t)-X_{1}^u(t))+i \sqrt{2} \beta (\bar{Y}_{j}^u(t)-\bar{Y}_{1}^u(t)  ) } \right ) \right]  \\
\nonumber & = \E \left[ F \left( \sum_{u \in \mathcal{H}_l(t)}  t^{\frac{3 \gamma}{2}} e^{-\gamma X_1^u(t)+i \sqrt{2} \beta \bar{Y}_1^u(t)   }  \sum_{j=1, \: t-\tau^u_{j,1}(t)<b}^{N^u(t)}\1_{\{X_j^u(t) \leq \frac{3}{2}\ln t +k \}}  e^{-\gamma (X_{j}^u(t)-X_{1}^u(t))+i \sqrt{2} \beta (\bar{Y}_{j}^u(t)-\bar{Y}_{1}^u(t)  ) } \right ) \right] \\
\label{ABORNE}&+ B_{t,l,k,b}
\end{align}
 where the remainder term $B_{t,l,k,b}$ is such that $|B_{t,l,k,b} |\leq ||F||_{\infty} \P(\mathcal{B}_{t,l,k,b})$ where $\mathcal{B}_{t,l,k,b}$ is defined by
\begin{equation*}
\mathcal{B}_{t,l,k,b} = \lbrace \exists j \in [|1, N^u(t)  |]; \: \tau^u_{j,1}(t) \leq t-b \: \text{and} \: X_j^u(t) \leq \frac{3}{2}\ln t +k    \rbrace.
\end{equation*}

Now for all $k' \geq k$, we have by lemma \ref{complextensionfinale} that 
\begin{align}
\nonumber&  \E \left[ F \left( \sum_{u \in \mathcal{H}_l(t)}  t^{\frac{3 \gamma}{2}} e^{-\gamma X_1^u(t)+i \sqrt{2} \beta \bar{Y}_1^u(t)   }  \sum_{j=1, \: t-\tau^u_{j,1}(t)<b}^{N^u(t)}\1_{\{X_j^u(t) \leq \frac{3}{2}\ln t +k \}}  e^{-\gamma (X_{j}^u(t)-X_{1}^u(t))+i \sqrt{2} \beta (\bar{Y}_{j}^u(t)-\bar{Y}_{1}^u(t)  ) } \right ) \right] =  \\
\nonumber &  \E \left[ F \left( \sum_{u \in \mathcal{H}_l(t)}  t^{\frac{3 \gamma}{2}} e^{-\gamma X_1^u(t)+i \sqrt{2} \beta \bar{Y}_1^u(t)   } \1_{\{X_1^u(t) \leq \frac{3}{2}\ln t +k \}} \sum_{j=1, \: t-\tau^u_{j,1}(t)<b}^{N^u(t)} \1_{\{X_j^u(t) \leq \frac{3}{2}\ln t +k' \}}  e^{-\gamma (X_{j}^u(t)-X_{1}^u(t))+i \sqrt{2} \beta (\bar{Y}_{j}^u(t)-\bar{Y}_{1}^u(t)  ) } \right ) \right] \\
\label{deuxsansdelta}&+ C_{t,l,k,k',b}
\end{align}
where $C_{t,l,k,k',b}$ is such that $\underset{t \to \infty}{\overline{\lim}}|C_{t,l,k,k',b}| \leq D_k$ where $D_k$ goes to $0$ when $k$ goes to infinity.

Now, in order to describe the limit, we need to introduce some notations. We consider $\mathcal{H}_l$ as a subset of $\N$. We introduce an i.i.d. sequence $(U_{u})_{u \in \N}$ of random variables uniformly distributed on the unit circle and an i.i.d. sequence $(B^{(u)})_{ u \in \N}$ of standard Brownian motions. We also consider an i.i.d. sequence $(\Gamma^{(u)})_{u \in \N}$ distributed like the backward path $Y$ of \cite{ABBS} and the associated Poisson jumps $((\tau_j^{(u)})_j)_{u \in N}$. Finally, given $u$ and conditionally to $(\Gamma^{(u)},(\tau_j^{(u)})_j)$, we consider an independent sequence of Point processes $\bar{\mathcal{N}}^{(u,j)}_{\Gamma^{(u)}(\tau_j^{(i)}),\tau_j^{(u)}}(dX,dY)$ of distribution that of $\bar{\mathcal{N}}_{x,s}(dX,dY) $ where $x=\Gamma^{(u)}(\tau_j^{(u)})$ and $s=\tau_j^{(u)}$.          

Now, using the convergence results of \cite{ABBS}, the $\E_{t,l,k,k',b}[ \cdots   ]$ term in the right hand side of  \eqref{deuxsansdelta} satisfies the following convergence

\begin{equation}\label{limitk}
 \underset{b \to \infty}{\lim} \underset{t \to \infty}{\lim} \E_{t,l,k,k',b}[ \cdots   ]
= \E \left [  F \left ( \sum_{u \in \mathcal{H}_l} \1_{\{W_{u}+l \leq k\}} \: e^{- \gamma (W_{u}+l)}  U_{u} Z_{W_u+l,k'}^{(u)}   \right )  \right ]
\end{equation}
 where $(W_u)_{u \in \mathcal{H}_l}$ is an i.i.d. sequence distributed like the asymptotic minimum of the (shifted) BBM, i.e. of distribution $W$ in \eqref{convBram}, and $Z_{W_u+l,k'}^{(u)}$ is an i.i.d. sequence given by 
\begin{equation*}
Z_{W_u+l,k'}^{(u)}= 1+\sum_{j \geq 1} e^{-\gamma \Gamma^{(u)}(\tau_j^{(u)}) - i \sqrt{2} \beta B^{(u)}(\tau_j^{(u)})}   < \1_{\{X+W_u +l\leq k'\}} e^{-\gamma X+i  \sqrt{2} \beta Y},  \bar{\mathcal{N}}^{(u,j)}_{\Gamma^{(u)}(\tau_j^{(u)}),\tau_j^{(u)}}(dX,dY) >
\end{equation*} 
The point that does not come out of the results of \cite{ABBS} is the appearance of the sequence $(U_{u})_{u \in \mathcal{H}_l}$ and the sequence $(B^{(u)}(\tau_j^{(u)})_j)_u$.  
Observe that if $\mathcal{G}$ is a standard Gaussian variable then $e^{i\alpha \mathcal{G}}$ converges in law as $\alpha \to \infty$ to a random variable uniformly distributed on the unit circle. We extend this elementary result to the following lemma:
\begin{lemma}
Assume that $(X^n)_n$ is a sequence of centered $\R^d$-valued Gaussian random vectors such that
\begin{equation}\label{fourier}
\forall x\in \R^d\setminus\{0\},\quad \E[|<x,X^n>|^2] \to \infty \text{ as }n\to\infty.
\end{equation}
Then the following convergence holds in law as $n\to\infty$
$$(e^{iX^n_1},\cdots,e^{iX^n_d})\to (U_1,\dots,U_d)$$ where $U_1,\dots,U_d$ are independent random variables uniformly distributed on the unit circle. 
\end{lemma}

\noindent {\it Proof.} Let us consider $d$   smooth functions $F_1,\dots F_d$ on the unit circle. We can write the Fourier expansion of the product $F_1\times F_d$ 
$$\forall x_1,\dots,x_d\in\R,\quad F_1(e^{ix_1})\times F_d(e^{ix_d})=\sum_{p\in\Z^d}c_p e^{i<p,x>}.$$
The sum is absolutely converging. We deduce
\begin{align*}
\E[F_1(e^{iX^n_1})\times F_d(e^{iX^n_d})]=\sum_{p\in\Z^d}c_p \E[e^{i<p,X^n>}]
\end{align*}
The relation $\E[e^{i<p,X^n>}]=e^{-\frac{1}{2}\E[|<p,X^n>|^2]}$ and assumption \eqref{fourier} imply that each term in the above sum, except for $p=0$, converges towards $0$. The dominated convergence theorem then entails that
$$\E[F_1(e^{iX^n_1})\times F_d(e^{iX^n_d})]\to c_0.$$
The result follows.\qed

%

Since $l$ is fixed here, conditionally on the Poisson process $N$ and the particles $X_j^u$, the sequence $(\bar{Y}^u_{1}(t))_{u}$ satisfies the assumptions of the above lemma as they are Brownian motions, the increments of which become independent after some time $t_\star$. Hence, since $\beta>0$, the sequence $(e^{i \sqrt{2} \beta \bar{Y}^u_{1}(t)})_u   $ converges in law as $t \to \infty$ to an i.i.d. sequence of random variables uniformly distributed on the unit circle and independent from all the other variables.

The other point to adress is the appearance of the sequence $(B^{(u)}(\tau_j^{(u)})_j)_u$.  Recall the following lemma:
\begin{lemma}\label{lemmejusti}
Let $(B(t))_{t \geq 0}$ be a standard Brownian motion. For all fixed $0 \leq t_1< \ldots <t_d$ and $a>0$, we have the following convergence in law 
\begin{equation*}
( e^{i a B(t)},e^{-i a( B(t)-B(t-t_1))}, \cdots,e^{- i a( B(t)-B(t-t_d))}  ) \underset{t \to \infty}{\rightarrow} (U,e^{-ia \tilde{B}(t_1)},\cdots ,e^{-ia \tilde{B}(t_d)} )
\end{equation*}
where $U$ is uniformly distributed on the circle and $\tilde{B}$ is a standard Brownian motion independent from $U$.
\end{lemma}
When $u$ is fixed, each term $\bar{Y}_{j}^u(t)-\bar{Y}_{1}^u(t)$ is the sum of $-(\bar{Y}_{1}^u(t)-\bar{Y}_{1}^u(\tau^u_{j,1}(t))  )$ and an independent branching part. Hence, conditionally on the Poisson process $N$ and the particles $X_j^u$, we can apply a straightforward variant of lemma \ref{lemmejusti} in the limit \eqref{limitk} since $b$ is fixed before taking the limit $t \to \infty$.

Now, we wish to take the limit in $k'$ in \eqref{limitk}. By lemma \ref{exislimitbis} and because the set $\mathcal{H}_l$ is finite, we have
\begin{equation*}
 \underset{k' \to \infty}{\lim} \underset{b \to \infty}{\lim}   \underset{t \to \infty}{\lim} \E_{t,l,k,k',b}[ \cdots   ]= \E \left [  F \left ( \sum_{u \in \mathcal{H}_l} \1_{\{W_{u}+l \leq k\}} \: e^{- \gamma (W_{u}+l)}  U_{u} Z^{(u)}   \right )  \right ]
\end{equation*}
 where $Z^{(u)}$ is an i.i.d. sequence given by (see lemma \ref{exislimit})
\begin{equation*}
Z^{(u)}= 1+\sum_{j \geq 1} e^{-\gamma \Gamma^{(u)}(\tau_j^{(u)}) - i \sqrt{2} \beta B^{(u)}(\tau_j^{(u)})}   <  e^{-\gamma X+i  \sqrt{2} \beta Y},  \bar{\mathcal{N}}^{(u,j)}_{\Gamma^{(u)}(\tau_j^{(u)}),\tau_j^{(u)}}(dX,dY) >.
\end{equation*}

Now, we wish to take the limit as $l$ goes to infinity. By the results of \cite{ABBS}, we get that  
\begin{align*}
& \underset{l \to \infty}{\lim}  \underset{k' \to \infty}{\lim} \underset{b \to \infty}{\lim} \underset{t \to \infty}{\lim}  \E_{t,l,k,k',b}[ \cdots   ]  = \E \left [  F \left ( \sum_{u  \geq 1} \1_{\{\Delta_u \leq k\}} \: e^{- \gamma \Delta_{u} } U_{u} Z^{(u)}   \right )  \right ],
\end{align*}       
where $(\Delta_u)_{u \geq 1}$ is a Poisson Point Process of intensity $c M' e^x dx $ where $M'$ is the limit of the derivative martingale.

To sum things up, we have proven that
\begin{align}
& \underset{l \to \infty}{\overline{\lim}} \underset{k' \to \infty}{\overline{\lim}}\underset{b \to \infty}{\overline{\lim}}\underset{t \to \infty}{\overline{\lim}} \left | \E \left [ F \left ( t^{\frac{3 \gamma}{2}} \sum_{i=1}^{N(t)}e^{-\gamma X_{i}(t)+i \sqrt{2} \beta \bar{Y}_{i}(t) }  \right ) \right ]-\E \left [  F \left ( \sum_{u  \geq 1} \1_{\{\Delta_u \leq k\}} \: e^{- \gamma \Delta_{u} } U_{u} Z^{(u)}   \right )  \right ]  \right |  \nonumber \\
& \leq \underset{l \to \infty}{\overline{\lim}} \underset{k' \to \infty}{\overline{\lim}} \underset{b \to \infty}{\overline{\lim}}\underset{t \to \infty}{\overline{\lim}} \left ( C_k+\E [|A_{t,l,k}| \wedge 1]+ |B_{t,l,k,b} | + |C_{t,l,k,k',b}|  \right ) \nonumber \\
& \leq C_k +D_k\label{yeah}.
\end{align}   
In fact, the bounds that we have obtained along the proof hold uniformly with respect to the functions $F$ that are bounded by $1$ and $1$-Lipschitz. Let $\mathcal{F}$ denote the space of such functions. We have thus proved 
\begin{align}
&\underset{t \to \infty}{\overline{\lim}} \underset{F \in \mathcal{F}}{\sup}\left | \E \left [ F \left ( t^{\frac{3 \gamma}{2}} \sum_{i=1}^{N(t)}e^{-\gamma X_{i}(t)+i \sqrt{2} \beta \bar{Y}_{i}(t) }  \right ) \right ]-\E \left [  F \left ( \sum_{u  \geq 1} \1_{\{\Delta_u \leq k\}} \: e^{- \gamma \Delta_{u} } U_{u} Z^{(u)}   \right )  \right ]  \right |  \nonumber \\
& \leq C_k +D_k\label{yeahbis}.
\end{align}   
 Now, we conclude by using the following trick. Recall that the sequence $ (t^{\frac{3 \gamma}{2}} \sum_{i=1}^{N(t)}e^{-\gamma X_{i}(t)+i \sqrt{2} \beta \bar{Y}_{i}(t) } )_{t}$ is tight. Indeed, by lemma \ref{complextension}, it suffices to show that for all fixed $k$ the sequence 
 \begin{equation*} 
 (t^{\frac{3 \gamma}{2}} \sum_{i=1}^{N(t)} \1_{\{X_i(t) \leq \frac{3}{2}\ln t +k \}} e^{-\gamma X_{i}(t)+i \sqrt{2} \beta \bar{Y}_{i}(t) } )_{t}\end{equation*}
 is tight. But this results from the real case \cite{ABBS,ABK1,ABK2,ABK3}  and the bound
 \begin{equation*}
| t^{\frac{3 \gamma}{2}} \sum_{i=1}^{N(t)} \1_{\{X_i(t) \leq \frac{3}{2}\ln t +k \}} e^{-\gamma X_{i}(t)+i \sqrt{2} \beta \bar{Y}_{i}(t) }  |  \leq | \lbrace  i; \: X_i(t) \leq \frac{3}{2}\ln t +k \rbrace  | e^{- \gamma (X_1(t) - \frac{3}{2}\ln t)}  
 \end{equation*}

 Since the sequence $ (t^{\frac{3 \gamma}{2}} \sum_{i=1}^{N(t)}e^{-\gamma X_{i}(t)+i \sqrt{2} \beta \bar{Y}_{i}(t) } )_{t}$ is tight, we can find a sequence $(t_j)_{j \geq 1}$ going to infinity and such that it converges in law towards a random variable. From this subsequence, we can extract an increasing subsequence $(t_{j_k})_{j_k \geq 1}$ such that for all $k$, we have 

\begin{align*}
& \underset{F \in \mathcal{F}}{\sup}\left | \E \left [ F \left ( t_{j_k}^{\frac{3 \gamma}{2}} \sum_{i=1}^{N(t_{j_k})}e^{-\gamma X_{i}(t_{j_k})+i \sqrt{2} \beta \bar{Y}_{i}(t_{j_k}) }  \right ) \right ]-\E \left [  F \left ( \sum_{u  \geq 1} \1_{\{\Delta_u \leq k\}} \: e^{- \gamma \Delta_{u} } U_{u} Z^{(u)}   \right )  \right ]  \right |  \nonumber \\
& \leq C_k +D_k +\frac{1}{k^2}\label{yeahbis}.
\end{align*}   
Hence, we conclude that $ \sum_{u  \geq 1} \1_{\{\Delta_u \leq k\}} \: e^{- \gamma \Delta_{u} } U_{u} Z^{(u)}$ converges in law as $k$ goes to infinity. We would like to identify this law. Let $x \in \R^2$. We denote the scalar product by $<,>$ and $\E_X$ expectation with respect to a variable $X$. By isotropy of the uniform law on the unit circle, the random variables $(<U_{u} Z^{(u)},x>)_u$ have the same laws as  $(|x||U_{u}^1| |Z^{(u)}|\epsilon_u)_u$ where $U_u^1$ is the first component of $U_u$ and $(\epsilon_u)_u$ is an independent family of i.i.d random variables with law $\P(\epsilon_u=1)=1-\P(\epsilon_u=-1)=\frac{1}{2}$. Recalling that $(\Delta_{u})_u$ is a Poisson point process with intensity $cM'e^z\,dz$, we have
\begin{align*}
 \E[   e^{i \sum_{u} \1_{\{\Delta_u \leq k\}} e^{-\gamma \Delta_u } <x,U_uZ^{(u)} > }   ]  & =  \E_{M'}  \left [   e^{c M'   \E_{U,Z}\big[   \int_{ \{v \leq k\}} ( e^{i e^{-\gamma v}   <x,UZ>} -1 )e^{v} dv\big] }  \right ]   \\
&=\E_{M'}  \left [   e^{c M'   \E_{U,Z,\epsilon}\big[   \int_{ \{v \leq k\}} ( e^{i e^{-\gamma v} |x|  |U^1| |Z|\epsilon} -1 )e^{v} dv\big] }  \right ]   \\
& = \E_{M'} \left [   e^{c\frac{M'}{  2   }   \E_{U,Z}\big[   \int_{ \{v \leq k\}} ( e^{i e^{-\gamma v} |x|  |U^1| |Z|} +e^{-i e^{-\gamma v} |x|  |U^1| |Z|}-2 )e^{v} dv\big] }  \right ]   \\
& = \E_{M'} \left [   e^{- c\frac{M'}{   \gamma } |x|^{1/\gamma} \E_{UZ}  \big[ |U^1Z |^{1/\gamma}  \int_{\{w \geq  e^{-\gamma k}|x|  |U^1| |Z| \}}(1-\cos(w) ) \frac{du}{w^{1+ \frac{1}{\gamma}}} \big]}   \right ]  .
\end{align*}
Then, by the monotone convergence theorem, we have the following convergence
\begin{equation*}
\underset{k \to \infty}{\lim}\E_{U,Z}  \left [ |U^1Z|^{1/\gamma}   \int_{\{ w \geq e^{-\gamma k}|x||U^1Z| \}  }(1-\cos(w) ) \frac{dw}{w^{1+ \frac{1}{\gamma}}}  \right ]=c_\gamma \: \E_{U,Z}  [ |U^1Z|^{1/\gamma} ].
\end{equation*} 
It is important to observe that the expectation $\E_{U,Z}  [ |U^1Z|^{1/\gamma} ]$ is necessarily finite,  otherwise the family $ \sum_{u  \geq 1} \1_{\{\Delta_u \leq k\}} \: e^{- \gamma \Delta_{u} } U_{u} Z^{(u)}$ could not converge in law as $k\to\infty$. In conclusion, there exists some constant $c(\gamma,\beta)<\infty$ such that for all $x$
\begin{equation}\label{limfinale}
 \underset{k \to \infty}{\lim} \E[   e^{i \sum_{u} \1_{\{\Delta_u \leq k\}} e^{-\gamma \Delta_u } <x,U_uZ^{(u)} > }   ] = \E_{M'} \left [   e^{- c(\gamma,\beta) M'  |x|^{1/\gamma} }  \right ]   .
\end{equation}
 Now, inequality \eqref{yeahbis}  yields that $ t^{\frac{3 \gamma}{2}} \sum_{i=1}^{N(t)}e^{-\gamma X_{i}(t)+i \sqrt{2} \beta \bar{Y}_{i}(t) }$ also converges in law as $t$ goes to infinity to the variable whose Fourier transform is defined by the right hand side of \eqref{limfinale}.

\subsection{Technical lemmas}

\subsubsection*{Study of the BBM at a fixed time}

In this technical subsection, we do not suppose that the particles are ordered and we will identify the interval $[|1, N(t)|]$ with a random tree. In particular, given two particles $i,j$, we will denote $\tau_{i,j}$ their splitting time and set $n_{\tau_{i,j}}$ to be the node of the random tree where the splitting occurs. We start with the following lemma which we will need in the next subsection  

\begin{lemma}\label{lemmamoment}
Let $\beta >0$ and $\gamma \in ]\frac{1}{2},1]$ such that $\gamma +\beta >1$. Then there exists some constant $C>0$ such that
\begin{equation*}    
\sup_{t \geq 0}\E[ ( \sum_{i,j=1}^{N(t)}  e^{- \gamma X_i(t) - \gamma X_j(t) -2 \beta^2 (t-\tau_{i,j}(t)  )}    )^{1/ (2\gamma)}  |  ] \leq C
\end{equation*}

\end{lemma}

\proof
For simplicity, we suppose $t$ is an integer. We have
\begin{align*}
& \E\left(  \left [ \sum_{i,j=1}^{N(t)} e^{-2 \beta^2 (t-\tau_{i,j}(t))}  e^{ -\gamma X_i(t)  } \ee^{ -\gamma X_j(t) }\right  ] ^{ 1/ (2 \gamma) }  \right)   \\
& \leq  \E\left(  \left [ \sum_{l=1}^{t}   e^{-2 \beta^2 (t-l)}     \sum_{\tau \in [l,l+1]}\ee^{ -2\gamma X_{n_\tau}(\tau) } \sum_{i,j; \: \tau_{i,j}=\tau }  e^{ -\gamma( X_i(t)+X_j(t)-2 X_{n_\tau}(\tau))} \right  ] ^{1/ (2\gamma)}  \right)   \\
& \leq   \E  \left [ \sum_{l=1}^{t}   e^{- \frac{\beta^2}{\gamma} (t-l)}     \sum_{\tau \in [l,l+1]} e^{  {\green -}  X_{n_\tau}(\tau) } \left(\sum_{i,j; \: \tau_{i,j}=\tau }  e^{ -\gamma( X_i(t)+X_j(t)-2 X_{n_\tau}(\tau))} \right)^{1/ (2\gamma)}  \right ]    \\
\end{align*}
We introduce for any $l>0$,  $\sigma_1^{(l)}<\sigma_2^{(l)}<...$ the times of successive branching after $l$. We have
\begin{align*}
&   \E \left [ \sum_{\tau \in [l,l+1]}   e^{-  X_{n_\tau}(\tau) } \left(\sum_{i,j; \: \tau_{i,j}=\tau }  e^{ -\gamma( X_i(t)+X_j(t)-2 X_{n_\tau}(\tau))} \right)^{1/ (2\gamma)} \right]\\
&= \sum_{p\geq 0}\E \left [ e^{- X_{n_{\sigma_p^{(l)}}}(\sigma_p^{(l)})}  \1_{\{ \sigma_p^{(l)}\leq l+1\}}  \left(\sum_{i,j; \: \tau_{i,j}=\sigma_p^{(l)} }  e^{ -\gamma( X_i(t)+X_j(t)-2 X_{n_{\sigma_p^{(l)}}}(\sigma_p^{(l)}))} \right)^{1/ (2\gamma)} \right]     \\
& \leq  \sum_{p\geq 0}\E \left [ e^{ - X_{n_{\sigma_p^{(l)}}}(\sigma_p^{(l)})}  \1_{\{ \sigma_p^{(l)}\leq l+1\}}  \left( \E [ \sum_{i,j; \: \tau_{i,j}=\sigma_p^{(l)} }  e^{ -\gamma( X_i(t)+X_j(t)-2 X_{n_{\sigma_p^{(l)}}}(\sigma_p^{(l)}))}   | \sigma_p^{(l)} ]\right)^{1/ (2\gamma)} \right]     \\
& =  \sum_{p\geq 0}\E \left [ e^{ -X_{n_{\sigma_p^{(l)}}}(\sigma_p^{(l)})}  \1_{\{ \sigma_p^{(l)}\leq l+1\}} e^{\frac{(1-\gamma)^2(t-\sigma_p^{(l)})}{\gamma}} \right]     \\
\end{align*} 
Hence by using (\ref{critical}), we get
\begin{align*}
& \E  \left [ \sum_{l=1}^{t}   e^{- \frac{\beta^2}{\gamma} (t-l)}     \sum_{\tau \in [l,l+1]} e^{ {\green -}  X_{n_\tau}(\tau) } \left(\sum_{i,j; \: \tau_{i,j}=\tau }  e^{ -\gamma( X_i(t)+X_j(t)-2 X_{n_\tau}(\tau))} \right)^{1/ (2\gamma)}  \right ]    \\
& \leq   \E  \left [ \sum_{l=1}^{t}   e^{ \frac{(1-\gamma)^2-\beta^2}{\gamma} (t-l)}     \sum_{\tau \in [l,l+1]} e^{ {\green -}  X_{n_\tau}(\tau) }  \right ]   \\
 & \leq C  \sum_{l=1}^{t}   e^{ \frac{(1-\gamma)^2-\beta^2}{\gamma} (t-l)}  \\
& \leq C  \\
\end{align*}
since $(1-\gamma)^2-\beta^2<0 $.

\qed

Now, we state an intermediate lemma which we will need to prove the important lemma \ref{complextension}. First, we introduce a few notations we will use in the sequel. For $L \in \R$, set $I_t(L):=[\frac{3}{2} \ln t-L,\frac{3}{2} \ln t-L+1]$. 
For any $i \in [|1,N(t)|]$ and $x \geq t/2$, we denote by $s(i,x) \in [t/2,x]$ the real which realizes the infimum of the trajectory on $[t/2,x]$
\begin{equation*}
X_i(s(i,x))= \inf_{u \in [t/2,x]} X_i(u).
\end{equation*}   
Then for any $k_1,k,\, v\in \N,\, L\in \R$ we define $\mathcal{Z}^{k_1,k}(v,L)$ the subset defined by 
\begin{align*}
& i\in \mathcal{Z}^{k_1,k}(v,L) \\
& \Longleftrightarrow  \\
& \inf_{s \leq t}  X_i(s) \geq -k_1 ,\, X_i(t) \in I_t(-k)  ,\, s(i,t) \in [v,v+1],\,  X_i(s(i,t)) \in I_t(L)  \\
\end{align*}

\begin{lemma}\label{intermediaire}
Let $\beta >0$ and $\gamma \in ]\frac{1}{2},1]$ be such that $\beta >1-\gamma$. Let $\kappa>0$ be such that $\frac{1}{2}<\gamma \kappa <\frac{3}{4} \wedge \gamma$. There exist $C,\delta>0$ and $\alpha \in ]0,1[$ such that, if $a_L=e^{\alpha L}$, then for any $L_0\in \N$, $k_1\in \N$, $t\geq 1$
\begin{align}
\nonumber  & \P\left( \sum_{L=L_0+1}^{2\ln t}\left| \sum_{i=1}^{N(t)}  \1_{\{i \in  \underset{k \in \mathbb{Z}}{\bigcup} \: \underset{v\in \{ \lfloor t/2 \rfloor,...,t-\lfloor a_L \rfloor\}}{\bigcup}\mathcal{Z}^{k_1,k}(v,L)\}} e^{ -\gamma( X_i(t)-\frac{3}{2}\ln t) +i \sqrt{2}\beta \bar{Y}_i(t)}\right| \geq \epsilon \right)
\\
\label{?3?} & \leq \epsilon^{- \kappa} C ((1+k_1) e^{-\delta L_0}+e^{k_1} \sum_{L=L_0+1}^{2\ln t}e^{-\delta t})
\end{align}

\end{lemma}

\proof
In the proof, for simplicity, we will suppose that $t/2$ and $a_L$ are integers. We denote $\tau_{i,j}:= \tau_{i,j}(t)$ the time where two particles $X_i$ and $X_j$ have split. Let $\kappa>0$ such that $\frac{1}{2}<\gamma \kappa <\frac{3}{4}\wedge \gamma$. According to the Markov property, then the sub-additivity, the probability in (\ref{?3?}) is smaller than
\begin{align}
\nonumber & \sum_{L=L_0+1}^{2\ln t}\epsilon^{- \kappa}\E\left(  \left| \sum_{i=1}^{N(t)}  \1_{\{i \in  \underset{k \in \mathbb{Z}}{\bigcup} \: \underset{v\in \{t/2, ... ,t-a_L\}}{\bigcup}\mathcal{Z}^{k_1,k}(v,L)\}} e^{ -\gamma( X_i(t)-\frac{3}{2}\ln t) +i \sqrt{2}\beta \bar{Y}_i(t)}\right|^{\kappa}  \right)
\\
\label{retouJensen} &\leq \sum_{L=L_0+1}^{2\ln t}\epsilon^{- \kappa}\E\left(  \left| \sum_{i=1}^{N(t)}  \1_{\{i \in  \underset{k \in \mathbb{Z}}{\bigcup} \: \underset{v\in \{t/2, ... ,t-a_L\}}{\bigcup}\mathcal{Z}^{k_1,k}(v,L)\}} e^{ -\gamma( X_i(t)-\frac{3}{2}\ln t) +i \sqrt{2}\beta \bar{Y}_i(t)}\right|^{2\kappa}  \right)^\frac{1}{2}.
\end{align} 
Let us study for any $L\in [L_0+1,2\ln t]$ the expectations in the right hand side of (\ref{retouJensen}). We take the conditional expectation according to the real part of the BBM, then via the Jensen inequality we deduce that 
\begin{align}
\nonumber & \E\left(  \left| \sum_{i=1}^{N(t)}  \1_{\{i \in  \underset{k \in \mathbb{Z}}{\bigcup} \: \underset{v\in \{t/2, ... ,t-a_L\}}{\bigcup}\mathcal{Z}^{k_1,k}(v,L)\}} e^{ -\gamma( X_i(t)-\frac{3}{2}\ln t) +i \sqrt{2}\beta \bar{Y}_i(t)}\right|^{2\kappa}  \right) \\
\nonumber & \leq  \E\left(  \left [ \sum_{i,j=1}^{N(t)} e^{-2 \beta^2 (t-\tau_{i,j}(t))} \1_{\{i,j \in  \underset{k \in \mathbb{Z}}{\bigcup} \: \underset{v\in \{t/2,...,t-a_L\}}{\bigcup}\mathcal{Z}^{k_1,k}(v,L)\}} \ee^{ -\gamma( X_i(t)-\frac{3}{2}\ln t)} \ee^{ -\gamma( X_j(t)-\frac{3}{2}\ln t)}\right  ] ^{ \kappa}  \right)   
\\
\label{popopop}& \leq  \E\left(  \left [ \sum_{l=1}^{t}   e^{-2 \beta^2 (t-l)}    \sum_{v=t/2}^{t-a_L} \sum_{\tau \in [l,l+1]} e^{ -2\gamma( X_{n_\tau}(\tau)-\frac{3}{2}\ln t)} \sum_{i,j; \: \tau_{i,j}=\tau } \1_{\{i \in  \underset{k \in \mathbb{Z}}{\bigcup}\mathcal{Z}^{k_1,k}(v,L)\}} e^{ -\gamma( X_i(t)+X_j(t)-2 X_{n_\tau}(\tau))} \right  ] ^{\kappa}  \right),
\end{align}
where in the first inequality we have applied Jensen's inequality with $x \mapsto x^2$ and $\E[. | X]$, the conditional measure with $(X_i)_{1 \leq i \leq N(t)}$ fixed. 
By sub-additivity of $x \mapsto x^\kappa$, this is smaller than 
\begin{eqnarray*}
\E\left(  \sum_{l=1}^{t}   e^{-2 \kappa \beta^2 (t-l)}    \sum_{v=t/2}^{t-a_L} \sum_{\tau \in [l,l+1]}\ee^{ -2\kappa \gamma( X_{n_\tau}(\tau)-\frac{3}{2}\ln t)} \1_{\{n_\tau \in A(v,l)\}} \left [\sum_{i,j; \: \tau_{i,j}=\tau }  e^{ -\gamma( X_i(t)+X_j(t)-2 X_{n_\tau}(\tau))} \right  ] ^{\kappa}  \right)
\end{eqnarray*}
where $n_\tau \in A(v,l)$ means:
\begin{eqnarray*}
 &  \inf_{s\leq \tau}{X_{n\tau}}(s)\geq -k_1  &  \quad     \text{if   }\quad l+1\leq t/2 ,
\\
 &  \inf_{s\leq \tau}{X_{n_\tau}}(s) \geq -k_1,\,   \inf_{s\in [t/2,\tau]}{X_{n_\tau}}(s) \geq  \frac{3}{2} \ln t-L  &  \quad     \text{if   }\quad t/2<l+1\leq v+2
                 \\
  &  \inf_{s\leq \tau}{X_{n_\tau}}(s) \geq -k_1,\,  s({n_\tau},\tau)\in [v,v+1],\, \inf_{s\in [t/2,\tau]}{X_{n_\tau}}(s) \in I_t(L)  &       \quad      \text{if   }\quad  v+1< l,
\end{eqnarray*}
where $s({n_\tau},\tau)$ satisfies $X_{n_\tau}(s({n_\tau},\tau))= \inf_{u \in [t/2,\tau]} X_{n_\tau}(u)$.   
By introducing, as in the proof of lemma \ref{lemmamoment}, for any $l>0$,  the times $\sigma_1^{(l)}<\sigma_2^{(l)}<...$ of successive branching after $l$, one can use the branching property at these times and Jensen's inequality (with $x \mapsto x^\kappa$) to get
\begin{eqnarray}
\nonumber &&  \sum_{l=1}^{t}   e^{-2 \kappa \beta^2 (t-l)}    \sum_{v=t/2}^{t-a_L} \E\left( \sum_{\tau \in [l,l+1]}\ee^{ -2\kappa \gamma( X_{n_\tau}(\tau)-\frac{3}{2}\ln t)}\1_{\{n_\tau \in A(l,v)\}}  \E \left( \sum_{i,j; \: \tau_{i,j}=\tau } e^{ -\gamma( X_i(t)+X_j(t)-2 X_{n_\tau}(\tau))}  \right)^{\kappa}  \right)
\\
\nonumber &&=  \sum_{l=1}^{t}   e^{-2 \kappa \beta^2 (t-l)}    \sum_{v=t/2}^{t-a_L} \E\left( \sum_{\tau \in [l,l+1]} e^{ -2\kappa \gamma( X_{n_\tau}(\tau)-\frac{3}{2}\ln t)}\1_{\{n_\tau \in A(v,l)\}}    \E\left( \sum_{i=1}^{N(t-\tau)} e^{-\gamma X_i(t-\tau)}   |  \tau \right)^{2\kappa} \right)
\\
\nonumber &&\leq C  \sum_{l=1}^{t}   e^{- \kappa \theta(\gamma,\beta) (t-l)}    \sum_{v=t/2}^{t-a_L} \E\left( \sum_{\tau \in [l,l+1]} e^{ -2\kappa \gamma( X_{n_\tau}(\tau)-\frac{3}{2}\ln t)}\1_{\{n_\tau \in A(v,l)\}}    \right)
\end{eqnarray}
where, in the last inequality, we have used the many-to-one lemma to evaluate $ \E\left( \sum_{i=1}^{N(t-\tau)} e^{-\gamma X_i(t-\tau)}  | \tau \right)$ and with $\theta(\gamma,\beta):= 2 (\beta^2-(1-\gamma)^2)>0$. Let us estimate $ \E\left( \sum_{\tau \in [l,l+1]} e^{ -2\kappa \gamma( X_{n_\tau}(\tau)-\frac{3}{2}\ln t)}\1_{\{n_\tau \in A(v,l)\}}    \right)$ according to the value of $v$ and $l$. For any $i\in \{1,...,N(l)\}$ we denote by $\Upsilon^{(i)}$ the set of all the branching times occurring along the BBM starting from $X_i(l)$.
\begin{align}
 \nonumber & \E\left( \sum_{\tau \in [l,l+1]} e^{ -2\kappa \gamma( X_{n_\tau}(\tau)-\frac{3}{2}\ln t)}\1_{\{n_\tau \in A(v,l)\}}    \right) \\
 \nonumber &= \E\left(\sum_{i=1}^{N(l)} e^{ -2\kappa \gamma( X_{i}(l)-\frac{3}{2}\ln t)} \sum_{\tau\in \Upsilon^{(i)},\tau\leq l+1} e^{ -2\kappa \gamma( X_{n_{\tau}}(\tau)-X_i(l))}  \1_{\{n_\tau \in A(v,l)\}}\right)
\\
\label{goodamount} &\leq \E\left(\sum_{i=1}^{N(l)} e^{ -2\kappa \gamma( X_{i}(l)-\frac{3}{2}\ln t)} \1_{\{ i \in B(v,l)\}}\right) \E\left(\sum_{\tau\leq 1} e^{ -2\kappa \gamma( X_{n_{\tau}}(\tau))} \right)  ,
\end{align}
where $i\in B(v,l)$ means
\begin{eqnarray*}
 &  \inf_{s\leq l}{X_i}(s)\geq -k_1  &  \quad     \text{if   }\quad l+1\leq t/2 ,
\\
 &  \inf_{s\leq l}{X_i}(s) \geq -k_1,\,   \inf_{s\in [t/2,l]}{X_i}(s) \geq  \frac{3}{2} \ln t-L  &  \quad     \text{if   }\quad t/2<l+1\leq v+2
                 \\
  &  \inf_{s\leq l}{X_i}(s) \geq -k_1,\,  s(i,l)\in [v,v+1],\, \inf_{s\in [t/2,l]}{X_i}(s) \in I_t(L)&       \quad      \text{if   }\quad  v+1< l.
\end{eqnarray*}
where $s(i,l)$ satisfies $X_{i}(s(i,l))= \inf_{u \in [t/2,l]} X_{i}(u)$.
We bound $\E\left(\sum_{\tau\leq 1} e^{ -2\kappa \gamma( X_{n_{\tau}}(\tau))} \right) $ by $C$ and we deduce by the many-to-one lemma and the Girsanov lemma that
\begin{align}
  \nonumber \E\left(\sum_{i=1}^{N(l)} e^{ -2\kappa \gamma( X_{i}(l)-\frac{3}{2}\ln t)} \1_{\{ i \in B(v,l)\}}\right) 
&\leq Ct^{3\kappa\gamma} e^l\E\left( \ee^{-2\kappa\gamma(\sqrt{2}B_l+2l)}\1_{\{\sqrt{2}B_\cdot+ 2\cdot \in B(v,l)\}}\right) \nonumber
\\
& \label{moodforl}=Ct^{3\kappa\gamma}\E\left(  e^{\sqrt{2}(1-2\kappa\gamma )B_l}\1_{\{\sqrt{2}B_\cdot\in B(v,l)\}} \right).\quad
\end{align}
where, in a slight abuse of notation, the condition $\sqrt{2} B_\cdot \in B(v,l)$ means that the trajectory satisfies the same conditions as $X_i$ when $i \in B(v,l)$.
Finally, we have established the bound
\begin{align}
\nonumber & \E\left(  \left| \sum_{i=1}^{N(t)}  \1_{\{i \in  \underset{k \in \mathbb{Z}}{\bigcup} \: \underset{v\in \{t/2, ... ,t-a_L\}}{\bigcup}\mathcal{Z}^{k_1,k}(v,L)\}} e^{ -\gamma( X_i(t)-\frac{3}{2}\ln t) +i \sqrt{2}\beta \bar{Y}_i(t)}\right|^{2\kappa}  \right) \\ 
& \leq  C  \sum_{l=1}^{t}   e^{- \kappa \theta(\gamma,\beta) (t-l)}    \sum_{v=t/2}^{t-a_L}    t^{3\kappa\gamma}\E\left(  e^{\sqrt{2}(1-2\kappa\gamma )B_l}\1_{\{\sqrt{2}B_\cdot\in B(v,l)\}} \right).   \label{bpose}
\end{align}
 Recall that $1-2\kappa\gamma<0$. According to the definition of $B(v,l)$, we divide the estimation of $\E\left(  e^{\sqrt{2}(1-2\kappa\gamma )B_l}\1_{\{\sqrt{2}B_\cdot\in B(v,l)\}} \right)$ in the following cases:

{-First case, $l+1\leq 3t/4$.}
\begin{equation}
\label{casou1}  t^{3\kappa\gamma}\E\left(  e^{\sqrt{2}(1-2\kappa\gamma )B_l}\1_{\{\sqrt{2}B_\cdot\in B(v,l)\}} \right)\leq ct^{3\kappa\gamma} \ee^{(1-2\kappa\gamma)k_1}.
\end{equation}

{-Second case, $3t/4\leq l\leq v+1$.}
\begin{align}
 &  t^{3\kappa\gamma} \E\left(  e^{\sqrt{2}(1-2\kappa\gamma )B_l}\1_{\{\sqrt{2}B_\cdot\in B(v,l)\}} \right)  \nonumber \\
&\leq  C e^{-(1-2\kappa\gamma)L}  \sum_{j\geq0} \ee^{\sqrt{2}(1-2\kappa\gamma)j} t^{\frac{3}{2}}\P\left(   \inf_{s\leq l}\sqrt{2}  B_s \geq -k_1,\,   \inf_{s\in [t/2,l]} \sqrt{2}B_s \geq  \frac{3}{2} \ln t-L,\,\sqrt{2} B_l\in I_t(L-j) \right)  \nonumber
 \\
 &\leq C e^{- (1-2\kappa\gamma)L}  \sum_{j\geq0} \ee^{ (1-2\kappa\gamma)j} t^{\frac{3}{2}} \frac{(1+k_1)(1+j)}{ t^{\frac{3}{2}}} \nonumber  \\
& \leq C (1+k_1) e^{- (1-2\kappa\gamma)L},  \label{casou2}
\end{align}
where we have used standard estimates for Brownian motion (see for example Lemmas 2.2 and 2.4 in \cite{AShi12}). 

{-Third case, $3t/4\leq v+1\leq l$}. By introducing $\sigma:= \inf\{s\geq t/2,\, B_s\leq \frac{3}{2} \ln t-L+1\}$, via the Markov property at time $\sigma$, we have for any $j\geq 0$
\begin{align*}
& \P\left(   \inf_{s\leq l}\sqrt{2} B_s \geq -k_1,\,  \sigma \in [v,v+1],\, \inf_{s\in [t/2,l]}\sqrt{2} B_s\in I_t(L),\, \sqrt{2}B_l\in I_t(L-j) \right)
\\
&\leq \E\left( 1_{\{  \inf_{s\leq v}\sqrt{2} B_s \geq -k_1,\,  \sigma \in [v,v+1]\}} \P_{B_\sigma-\frac{3}{2} \ln t+L}\left( \inf_{s\in[t/2-\sigma,l-\sigma]} \sqrt{2}B_s\geq -1,\, \sqrt{2}B_{l-\sigma}\in [j,j+1] \right)\right)
\\
&\leq C \frac{1+j}{(l-v)^\frac{3}{2}}\P\left( \inf_{s\leq v}\sqrt{2} B_s \geq -k_1,\,  \sigma \in [v,v+1]\right) \\
& \leq C \frac{(1+j)(1+k_1)}{(l-v)^\frac{3}{2} t^\frac{3}{2}},
\end{align*}
where we have used standard estimates for Brownian motion (see for example Lemmas 2.2 and 2.4 in \cite{AShi12}). Thus we deduce that
\begin{align}
 \nonumber t^{3\kappa\gamma} \E\left(  e^{\sqrt{2}(1-2\kappa\gamma )B_l}\1_{\{\sqrt{2}B_\cdot\in B(v,l)\}} \right)  &  \leq C  e^{-(1-2\kappa\gamma)L}  \sum_{j\geq0} e^{\sqrt{2}(1-2\kappa\gamma)j} t^{\frac{3}{2}}  \frac{(1+j)(1+k_1)}{(l-v)^\frac{3}{2} t^\frac{3}{2}} \\ & \leq   C \frac{(1+k_1)}{(l-v)^\frac{3}{2}} e^{- (1-2\kappa\gamma)L}.  \label{casou3} 
\end{align}
Going back to (\ref{bpose}), and by combining (\ref{casou1}), (\ref{casou2}) and (\ref{casou3}) we get
\begin{align}
\nonumber & \E\left(  \left| \sum_{i=1}^{N(t)}  \1_{\{i \in  \underset{k \in \mathbb{Z}}{\bigcup} \: \underset{v\in \{t/2, ... ,t-a_L\}}{\bigcup}\mathcal{Z}^{k_1,k}(v,L)\}} e^{ -\gamma( X_i(t)-\frac{3}{2}\ln t) +i \sqrt{2}\beta \bar{Y}_i(t)}\right|^{2\kappa}  \right)  \\
 & \leq  C   \sum_{l=1}^{t}   e^{- \kappa \theta(\gamma,\beta) (t-l)}    \sum_{v=t/2}^{t-a_L}    t^{3\kappa\gamma}\E\left(  e^{\sqrt{2}(1-2\kappa\gamma )B_l}\1_{\{\sqrt{2}B_\cdot\in B(v,l)\}} \right) \nonumber \\  
& \leq (A)+(B)+(C)+(D)  \label{goback}
\end{align}
where 
\begin{align}
\nonumber (A)&:=\sum_{l=1}^{3t/4} e^{- \kappa \theta(\gamma,\beta) (t-l)} \sum_{v=t/2}^{t-a_L} t^{3\kappa\gamma} \E\left(  e^{\sqrt{2}(1-2\kappa\gamma )B_l}\1_{\{\sqrt{2}B_\cdot\in B(v,l)\}} \right)
\\
\nonumber &\leq C \sum_{l=1}^{3t/4} e^{- \kappa \theta(\gamma,\beta) (t-l)} \sum_{v=t/2}^{t-a_L}  t^{3\kappa\gamma} e^{(1-2\kappa\gamma)k_1}  \\
& \label{**1} \leq C t^{3\kappa\gamma +1} e^{k_1} e^{-\kappa\theta(\gamma,\beta)t/4}
\end{align}
and
\begin{align}
\nonumber (B)&:=    \sum_{l=3t/4}^{t-\frac{a_L}{2}}e^{ -\kappa\theta(\gamma,\beta)(t-l)}  \underset{v=l+1}{\overset{t-a_L}{\sum}}   t^{3\kappa\gamma} \E\left(  \ee^{\sqrt{2}(1-2\kappa\gamma )B_l}\1_{\{\sqrt{2}B_\cdot\in B(v,l)\}} \right)
\\
\nonumber  &\leq C \sum_{l=3t/4}^{t-\frac{a_L}{2}} e^{ -\kappa\theta(\gamma,\beta)(t-l)}  (t-l)  (1+k_1) e^{- (1-2\kappa\gamma)L} \\
\label{**2}  &\leq C(1+k_1) e^{-\kappa\theta(\gamma,\beta)a_L/4-(1-2\kappa\gamma)L},
\end{align}
where we have used the inequality $x e^{-x}\leq C e^{-x/2}$ for any $x\geq 1$,
\begin{align}
\nonumber (C)&:=  \sum_{l=3t/4}^{t-\frac{a_L}{2}} e^{ -\kappa\theta(\gamma,\beta)(t-l)}  \underset{v=t/2}{\overset{(l-1) \wedge (t-a_L)}{\sum}}t^{3\kappa\gamma} \E\left(  e^{\sqrt{2}(1-2\kappa\gamma )B_l}\1_{\{\sqrt{2}B_\cdot\in B(v,l)\}} \right)
\\
\nonumber &\leq C \sum_{l=3t/4}^{t-\frac{a_L}{2}} e^{ -\kappa\theta(\gamma,\beta)(t-l)}  \underset{v=t/2}{\overset{(l-1) \wedge (t-a_L)}{\sum}}  \frac{(1+k_1)}{(l-v)^\frac{3}{2}} e^{- (1-2\kappa\gamma)L}
\\
\nonumber &\leq C \sum_{l=3t/4}^{t-\frac{a_L}{2}}e^{ -\kappa\theta(\gamma,\beta)(t-l)}     {(1+k_1)} e^{- (1-2\kappa\gamma)L}  \\
& \label{**3} \leq C(1+k_1) e^{-\kappa\theta(\gamma,\beta)a_L/2-(1-2\kappa\gamma)L}  
\end{align}
\begin{align}
\nonumber (D)&:=   \sum_{l=t-\frac{a_L}{2}}^{t} e^{ \kappa\theta(\gamma,\beta)(t-l)}  \underset{v=t/2}{\overset{t-a_L}{\sum}}  t^{3\kappa\gamma} \E\left(  e^{\sqrt{2}(1-2\kappa\gamma )B_l}\1_{\{\sqrt{2}B_\cdot\in B(v,l)\}} \right)
\\
\nonumber &\leq C   \sum_{l=t-\frac{a_L}{2}}^{t} e^{ \kappa\theta(\gamma,\beta)(t-l)}  \underset{v=t/2}{\overset{t-a_L}{\sum}}   \frac{(1+k_1)}{(l-v)^\frac{3}{2}} e^{- (1-2\kappa\gamma)L}
\\
\nonumber &\leq  C \sum_{l=t-\frac{a_L}{2}}^{t}e^{ -\kappa\theta(\gamma,\beta)(t-l)}     {(1+k_1)} \frac{e^{- (1-2\kappa\gamma)L}}{a_L^\frac{1}{2}}  \\
\label{**4} & \leq C (1+k_1) a_L^{-\frac{1}{2}}e^{-(1-2\kappa\gamma)L},  
\end{align}
where we have used the inequality $\sum_{i\geq j}\frac{1}{i^\frac{3}{2}}\leq C i^{-\frac{1}{2}}$ for any $j \geq 1$.

Recall here that the condition on $\kappa$ ensures that $2\kappa\gamma-1>0$ and $\frac{1}{2}-(2\kappa\gamma -1)>0$. Hence we can find $\alpha \in ]0,1[$ and $\delta>0$ such that $\frac{\alpha}{2}-(2\kappa\gamma -1)>\delta$ leading to $a_L^{-\frac{1}{2}}e^{-(1-2\kappa\gamma)L} \leq e^{- \delta L}$. We also suppose that $\delta< \theta(\gamma,\beta)/4$. Then combining (\ref{goback}) with (\ref{**1}), (\ref{**2}), (\ref{**3}), (\ref{**4}),  it is plain to deduce that for any $t,\, L>0$
\begin{align}
\nonumber & \E\left(  \left| \sum_{i=1}^{N(t)}  \1_{\{i \in  \underset{k \in \mathbb{Z}}{\bigcup} \: \underset{v\in \{t/2, ... ,t-a_L\}}{\bigcup}\mathcal{Z}^{k_1,k}(v,L)\}} e^{ -\gamma( X_i(t)-\frac{3}{2}\ln t) +i \sqrt{2}\beta \bar{Y}_i(t)}\right|^{2\kappa}  \right)  \\
\nonumber & \leq C  \left(  t^{3\kappa\gamma +1} e^{k_1}e^{-\kappa\theta(\gamma,\beta)t/4}+(1+k_1)e^{- (1-2\kappa\gamma)L} \left( e^{-\kappa\theta(\gamma,\beta)a_L/4 }   +e^{-\kappa\theta(\gamma,\beta)a_L/2}+  a_L^{-\frac{1}{2}} \right) \right)  \\
\label{retouback22} & \leq C  ( (1+k_1)e^{-  \delta L} +e^{k_1}e^{-\delta t}).
\end{align}
  Going back to (\ref{retouJensen}), by using (\ref{retouback22}) for any $L\in [L_0+1,2\ln t]$ we get:
\begin{align*}
&  \P\left( \sum_{L=L_0+1}^{2\ln t}\left| \sum_{i=1}^{N(t)}  \1_{\{i \in  \underset{k \in \mathbb{Z}}{\bigcup} \: \underset{v\in \{ \lfloor t/2 \rfloor,...,t-\lfloor a_L \rfloor\}}{\bigcup}\mathcal{Z}^{k_1,k}(v,L)\}} e^{ -\gamma( X_i(t)-\frac{3}{2}\ln t) +i \sqrt{2}\beta \bar{Y}_i(t)}\right| \geq \epsilon \right) \\
&\leq \sum_{L=L_0+1}^{2\ln t}  \epsilon^{-\kappa}\sqrt{C ( (1+k_1)e^{-  \delta L} +e^{k_1}e^{-\delta t})}
\\
&\leq \epsilon^{- \kappa} C (\sqrt{(1+k_1)} e^{-\delta L_0/2}+e^{k_1/2} \sum_{L=L_0+1}^{2\ln t}e^{-\delta t/2})
\end{align*}
which is the desired result.

\qed

Here we prove the main lemma of this subsection

\begin{lemma}\label{complextension}
Let $\beta >0$ and $\gamma >\frac{1}{2}$ such that $\beta >(1-\gamma)_+$. Then we have
for any $\epsilon>0$, 
\begin{equation}
\label{limlim} \underset{k\to\infty}{\lim}\,  \underset{t \to \infty}{\overline{\lim}}\,  \P\left(\left| \sum_{i=1}^{N(t)}  \1_{\{X_i(t)-\frac{3}{2}\ln t \geq k \}} e^{ -\gamma( X_i(t)-\frac{3}{2}\ln t) +i \sqrt{2}\beta \bar{Y}_i(t)}\right| \geq \epsilon \right) \to 0
\end{equation}
\end{lemma}

\noindent{\it Proof of Lemma \ref{complextension}.}
If $\gamma>1$,we have the following obvious bound
\begin{equation*}
\left| \sum_{i=1}^{N(t)}  \1_{\{X_i(t)-\frac{3}{2}\ln t \geq k \}} e^{ -\gamma( X_i(t)-\frac{3}{2}\ln t) +i \sqrt{2}\beta \bar{Y}_i(t)}\right| \leq  \sum_{i=1}^{N(t)}  \1_{\{X_i(t)-\frac{3}{2}\ln t \geq k \}} e^{ -\gamma( X_i(t)-\frac{3}{2}\ln t)},
\end{equation*}
and therefore it suffices to adapt (from the branching random walk to the BBM) the proof of Proposition 4.6 in \cite{madaule} to obtain the result.

Hence, in the sequel, we suppose that $\gamma \in ]\frac{1}{2},1]$. For simplicity, we suppose $t/2$ is an integer. Recall from a minor adaptation of \cite{aidekon} that for any $k_1\geq 0$, we have
\begin{equation}
\label{minminmin}
\P(   \inf_{t \geq 0}  \inf_{1 \leq i \leq N(t)} X_i(t) \leq -k_1  )  \leq e^{-k_1}.
\end{equation}
For any $k_1\geq 0$, $L\geq 0$, the probability in \eqref{limlim} is less or equal than  
\begin{equation}
 e^{-k_1}+ \P\left(\left| \sum_{i=1}^{N(t)}  \1_{\{  \inf_{s \leq t}  X_i(s) \geq -k_1 , X_i(t)-\frac{3}{2}\ln t \geq k \}} e^{ -\gamma( X_i(t)-\frac{3}{2}\ln t) +i\sqrt{2} \beta \bar{Y}_i(t)}\right| \geq \epsilon \right)   
\end{equation}
Clearly, we have
\begin{equation}
\label{consAID} \{ i \in [|1,N(t)|],\,  \inf_{s \leq t}  X_i(s) \geq -k_1 , X_i(t)-\frac{3}{2}\ln t \in [k,k+1]  \}= \underset{v\in \{t/2,...,t\}}{\bigcup}\underset{L\in \{-k,...,2\ln t\}}{\bigcup}\mathcal{Z}^{k_1,k}(v,L) 
\end{equation}

We consider $\delta>0$ and $\alpha \in ]0,1[$ according to lemma \ref{intermediaire} (recall that $a_L=e^{\alpha L}$ that we will also suppose to be an integer for simplicity). According to the proof of Lemma 3.3 in \cite{aidekon} (or more precisely its analogue to BBM), we know that for any $t,\, k_1, L\in \N^*$, 
\begin{equation}
\label{AidadaptBBM}\P\big( \exists i \in [|1,N(t)|],\, i \in  \underset{k \in \mathbb{Z}}{\bigcup} \: \underset{v\in \{ t-a_L,...,t\}}{\bigcup}\mathcal{Z}^{k_1,k}(v,L)\big)\leq C (1+k_1) a_L e^{-L}.
\end{equation}
This inequality is useful for $L$ large.

Now, according to lemma \ref{intermediaire}, we have for any $k_1$, $t\geq 1$ and $L>0$,
\begin{align}
 & \P\left(\left| \sum_{i=1}^{N(t)}  \1_{\{X_i(t)-\frac{3}{2}\ln t \geq k \}} e^{ -\gamma( X_i(t)-\frac{3}{2}\ln t) +i \sqrt{2}\beta \bar{Y}_i(t)}\right| \geq \epsilon \right) \nonumber  \\
& \leq e^{-k_1} + \P\left(\left| \sum_{i=1}^{N(t)}  \1_{\{  \inf_{s \leq t}  X_i(s) \geq -k_1 , X_i(t)-\frac{3}{2}\ln t \geq k \}} e^{ -\gamma( X_i(t)-\frac{3}{2}\ln t) +i\beta\sqrt{2}\: \bar{Y}_i(t)}\right| \geq \epsilon \right)   \nonumber
\\
&\leq e^{-k_1}  + \P\left(\left| \sum_{i=1}^{N(t)}  \1_{\{   \inf_{s\leq t} X_i(s)\geq -k_1,\, \inf_{s\in [t/2,t]}X_i(s)\geq \frac{3}{2} \ln t-L,\, X_i(t)\geq \frac{3}{2} \ln t+k \}} e^{ -\gamma( X_i(t)-\frac{3}{2}\ln t) +i\beta\sqrt{2}\: \bar{Y}_i(t)}\right| \geq \epsilon \right)  \nonumber
\\
&+ \sum_{L'=L+1}^{2\ln t}\P\left( \exists i \in [|1,N(t)|],\, i \in  \underset{k \in \mathbb{Z}}{\bigcup} \: \underset{v\in \{ t-a_L,...,t\}}{\bigcup}\mathcal{Z}^{k_1,k}(v,L')\right)  \nonumber
\\
&+\P\left(\sum_{L'=L+1}^{2\ln t}\left| \sum_{i=1}^{N(t)}  \1_{\{i \in  \underset{k \in \red{\mathbb{Z}}}{\bigcup} \: \underset{v\in \{ t/2,...,t-a_L\}}{\bigcup}\mathcal{Z}^{k_1,k}(v,L')\}} e^{ -\gamma( X_i(t)-\frac{3}{2}\ln t) +i\beta\sqrt{2}\: \bar{Y}_i(t)}\right| \geq \epsilon \right) \nonumber
\\
&\leq  e^{-k_1} + C(1+k_1) e^{-\delta L}  + C e^{k_1} \ln t  e^{-\delta t} \nonumber
\\
&+  \P\left(\left| \sum_{i=1}^{N(t)}  \1_{\{   \inf_{s\leq t} X_i(s)\geq -k_1,\, \inf_{s\in [t/2,t]}X_i(s)\geq \frac{3}{2} \ln t-L,\, X_i(t)\geq \frac{3}{2} \ln t+k \}} e^{ -\gamma( X_i(t)-\frac{3}{2}\ln t) +i\beta\sqrt{2}\: \bar{Y}_i(t)}\right| \geq \epsilon \right),  \label{tttt}
\end{align}
where, in the last inequality, we have used the bound $\sum_{L'=L+1}^\infty e^{-\delta L'} \leq C e^{-\delta L}$.

Thus in order to prove (\ref{limlim}), it remains to study for $t,\, L>0$, 
\begin{equation}
\label{stud}  \P\left(\left| \sum_{i=1}^{N(t)}  \1_{\{   \inf_{s\leq t} X_i(s)\geq -k_1,\, \inf_{s\in [t/2,t]}X_i(s)\geq \frac{3}{2} \ln t-L,\, X_i(t)\geq \frac{3}{2}\ln t+k \}} e^{ -\gamma( X_i(t)-\frac{3}{2}\ln t) +i\beta\sqrt{2}\: \bar{Y}_i(t)}\right| \geq \epsilon \right)   .
\end{equation}
According to the Markov inequality and Jensen's inequality, the probability in (\ref{stud})  is smaller than 
\begin{align}
&\nonumber \epsilon^{-2}\E (  \sum_{l=1}^{t}   e^{-2 \beta^2 (t-l)}    \sum_{\tau \in [l,l+1]} e^{ -2\gamma( X_{n_\tau}(\tau)-\frac{3}{2}\ln t)} \sum_{i,j; \: \tau_{i,j}=\tau } \1_{\{   \inf_{s\leq t} X_i(s)\geq -k_1,\, \inf_{s\in [t/2,t]}X_i(s)\geq \frac{3}{2}\ln t-L,\, X_i(t)\geq \frac{3}{2}\ln t+k \}}  
\\
&\label{proproprop} \times e^{ -\gamma( X_i(t)+X_j(t)-2 X_{n_\tau}(\tau))}  )
\\
&\nonumber  = \epsilon^{-2}\E\left(  \sum_{l=1}^{t-e^{L}} ...\right) +\epsilon^{-2}\E\left(  \sum_{l=t-e^L+1}^{t} ... \right) 
\end{align}
Again by introducing for any $l>0$,  $\sigma_1^{(l)}<\sigma_2^{(l)}<...$ the times of successive branching after $l$, by the branching property at these time we can write:
\begin{align}
& \nonumber \epsilon^{-2}\E\left(  \sum_{l=t-e^L+1}^{t} ... \right)  \\
& \leq C \epsilon^{-2}\E  (  t^{3 \gamma }\sum_{l=t-e^L+1}^{t}   e^{[ (1-\gamma)^2- 2\beta^2 ](t-l)}     \sum_{\tau \in [l,l+1]}  \sum_{i > \tau; } \1_{\{   \inf_{s\leq t} X_i(s)\geq -k_1,\, \inf_{s\in [t/2,t]}X_i(s)\geq \frac{3}{2}\ln t-L,\, X_i(t)\geq \frac{3}{2}\ln t+k \}}   \nonumber  \\
& \times  e^{ -\gamma( X_i(t) + X_{n_\tau}(\tau))} ),  \label{pourlasuite}
\end{align}
where $i> \tau$ means that $\tau$ is a spliting time of particle $i$. In the above equality, we have averaged out the trajectory of particle $j$ on the interval $[\tau,t]$. We also have 
\begin{equation}
\label{secood2}
\epsilon^{-2}\E\left(  \sum_{l=1}^{t-e^{L}} ...\right)  \leq   \epsilon^{-2}\E\left(  \sum_{l=1}^{t-e^L}   e^{- \theta(\beta,\gamma) (t-l)}    \sum_{\tau \in [l,l+1]} e^{ -2\gamma( X_{n_\tau}(\tau)-\frac{3}{2}\ln t)}\1_{\{n_\tau \in A(l) \}} \right)
\end{equation}
where we have averaged out the particles $i,j$ on $[\tau,t]$ and $n_\tau \in A(l)$ means:
\begin{eqnarray*}
 &  \inf_{s\leq \tau}{X_{n_\tau}}(s)\geq -k_1  &  \quad     \text{if   }\quad l+1\leq 3t/4 ,
\\
 &  \inf_{s\leq \tau}{X_{n_\tau}}(s) \geq -k_1,\,   \inf_{s\in [t/2,\tau]}{X_{n_\tau}}(s) \geq  \frac{3}{2}\ln t -L  &  \quad     \text{if   }\quad 3t/4<l+1\leq t- e^L
\end{eqnarray*}
First let us bound the term in (\ref{secood2}). By reasoning as in (\ref{goodamount}) and (\ref{moodforl}) we have
\begin{eqnarray*}
\epsilon^{-2}\E\left(  \sum_{l=1}^{t-e^{L}} ...\right)&\leq& C \sum_{l=1}^{t-e^{L}} e^{- \theta(\beta,\gamma) (t-l)}\E\left(\sum_{i=1}^{N(l)} e^{ -2 \gamma( X_{i}(l)-\frac{3}{2}\ln t)} \1_{\{ i \in A(l)\}}\right) \E\left(\sum_{\tau\leq 1} e^{ -2\gamma( X_{n_{\tau}}(\tau))} \right)
\\
&\leq& C \sum_{l=1}^{t-e^{L}} e^{- \theta(\beta,\gamma) (t-l)}t^{3 \gamma}\E\left(  e^{\sqrt{2}(1-2 \gamma )B_l}\1_{\{\sqrt{2}B_\cdot\in A(l)\}} \right)
\end{eqnarray*}
where $\sqrt{2}B_\cdot \in A(l)$ means
\begin{eqnarray*}
 &  \sqrt{2}\inf_{s\leq l}{B_s}\geq -k_1  &  \quad     \text{if   }\quad l+1\leq 3t/4 ,
\\
 &  \sqrt{2} \inf_{s\leq l}B_s  \geq -k_1,\,  \sqrt{2} \inf_{s\in [t/2,l]}B_s \geq  \frac{3}{2}\ln t-L  &  \quad     \text{if   }\quad 3t/4<l+1\leq t-\ee^L.
\end{eqnarray*}

Then it follows that
\begin{align*}
 \epsilon^{-2}\E\left(  \sum_{l=1}^{t-e^{L}} ...\right)  \leq & C \sum_{l=1}^{3t/4}
e^{-\theta(\beta,\gamma)(t-l)} t^{3\gamma} e^{(1-2\gamma)k_1}+ C
\sum_{l=3t/4}^{t- e^L} e^{-\theta(\beta,\gamma)(t-l)}    e^{-(1-2\gamma)L}
\times
\\
& \sum_{j\geq0} e^{\sqrt{2}(1-2 \gamma)j} t^{\frac{3}{2}}\P\left(  
\inf_{s\leq l}\sqrt{2}  B_s \geq -k_1,\,   \inf_{s\in [t/2,l]} \sqrt{2}B_s \geq 
\frac{3}{2} \ln t-L,\,\sqrt{2} B_l\in I_t(L-j) \right)
 \\
  \leq & Ct^{3\gamma} e^{(1-2\gamma) k_1} e^{-\theta(\beta,\gamma)t/4} +C 
e^{-(1-2\gamma)L-\theta(\beta,\gamma)e^L} (1+k_1).
\end{align*}

Now we need to bound the term \eqref{pourlasuite}. We can bound the term in \eqref{pourlasuite} by
\begin{align*}
& \epsilon^{-2}\sum_{l=t-e^L+1}^{t}   e^{[ (1-\gamma)^2- 2\beta^2 ](t-l)} t^{3\gamma}
\E\left( \sum_{i=1}^{N(t)}   \1_{\{   \inf_{s\leq t} X_i(s)\geq -k_1,\, \inf_{s\in
[t/2,t]}X_i(s)\geq \frac{3}{2} \ln t-L,\, X_i(t)\geq \frac{3}{2} \ln t+k \}}     e^{ -\gamma( X_i(t) +
\inf_{s\in [l,l+1]} X_{i}(s))}  \right)
\\
& = \epsilon^{-2}\sum_{l=t-e^L+1}^{t}  t^{3\gamma}  e^{( (1-\gamma)^2- 2\beta^2 )(t-l)}  \\
& \times 
\E\left( e^{\sqrt{2} B_t-\gamma(\sqrt{2} B_t +\inf_{s\in [l,l+1]} \sqrt{2} B_s)}   \1_{\{   \inf_{s\leq t}
\sqrt{2} B_s\geq -k_1,\, \inf_{s\in [t/2,t]}\sqrt{2} B_s\geq \frac{3}{2} \ln t -L,\, \sqrt{2} B_t\geq \frac{3}{2} \ln t+k \}}    \right)  \\
& \leq  \epsilon^{-2} t^{3/2} e^{(1-2\gamma)k}  \sum_{l=t-e^L+1}^{t}   e^{[ (1-\gamma)^2- 2\beta^2 ](t-l)}  \\
& \times\sum_{j=1}^\infty e^{(1-2 \gamma)j} \E\left( e^{-\gamma(\inf_{s\in [l,l+1]} \sqrt{2} B_s-\sqrt{2}B_t)}   \1_{\{   \inf_{s\leq t}
\sqrt{2} B_s\geq -k_1,\, \inf_{s\in [t/2,t]}\sqrt{2} B_s\geq \frac{3}{2} \ln t -L,\, \sqrt{2} B_t \in I_t(-j-k) \}}    \right)  \\
&  \leq C(L) \epsilon^{-2} t^{3/2} e^{(1-2\gamma)k}  \sum_{l=t-e^L+1}^{t}   e^{((1-\gamma)^2- 2\beta^2 )(t-l)}  \\
& \times 
\sum_{j=1}^\infty e^{(1-2 \gamma)j}  \frac{(1+k_1)(L+j+k)}{t^{3/2}}  \\ 
& \leq C(L) \epsilon^{-2} (1+k_1)(L+k) e^{(1-2\gamma)k} 
\end{align*}
where $C(L)$ is a constant depending on $L$ and we have used standard estimates on Brownian motion (see for example Lemmas 2.2 and 2.4 in \cite{AShi12}). In conclusion, we have the following bound
\begin{align}
& \P\left(\left| \sum_{i=1}^{N(t)}  \1_{\{   \inf_{s\leq t} X_i(s)\geq -k_1,\, \inf_{s\in [t/2,t]}X_i(s)\geq \frac{3}{2} \ln t-L,\, X_i(t)\geq \frac{3}{2} \ln t+k \}} e^{ -\gamma( X_i(t)-\frac{3}{2}\ln t) +i\beta\sqrt{2}\: \bar{Y}_i(t)}\right| \geq \epsilon \right)  \nonumber \\
& \leq Ct^{3\gamma} e^{(1-2\gamma) k_1} e^{-\theta(\beta,\gamma)t/4} +C 
e^{-(1-2\gamma)L-\theta(\beta,\gamma)e^L} (1+k_1)+ C(L) \epsilon^{-2} (1+k_1)(L+k) e^{(1-2\gamma)k}. \label{ttttt}
\end{align} 
Gathering \eqref{tttt} and \eqref{ttttt}, we finally obtain
\begin{align*}
& \P\left(\left| \sum_{i=1}^{N(t)}  \1_{\{X_i(t)-\frac{3}{2}\ln t \geq k \}} e^{ -\gamma( X_i(t)-\frac{3}{2}\ln t) +i \sqrt{2}\beta \bar{Y}_i(t)}\right| \geq \epsilon \right)  \\
& \leq  e^{-k_1} + C(1+k_1) e^{-\delta L}  + C e^{k_1} \ln t  e^{-\delta t} \\
& + Ct^{3\gamma} e^{(1-2\gamma) k_1} e^{-\theta(\beta,\gamma)t/4} +C   e^{-(1-2\gamma)L-\theta(\beta,\gamma)e^L} (1+k_1)+ C(L) \epsilon^{-2} (1+k_1)(L+k) e^{(1-2\gamma)k} .
\end{align*} 
Now, one concludes by letting $t \to \infty$ and then choosing successively $k_1,L,k$.

\qed

\begin{lemma}\label{complextensionfinale}
We have the following limit for all $\epsilon>0$
\begin{equation*}
\underset{k \to \infty}{\lim}\sup_{k' \geq k} \sup_{b,\,l \geq 0} \underset{t \to \infty}{\overline{\lim}}\P \left ( \left |   \sum_{u \in \mathcal{H}_l(t)}  t^{\frac{3 \gamma}{2}} \1_{\{X_1^u(t) \leq \frac{3}{2}\ln t +k \}}  \sum_{j=1, \: t-\tau^u_{j,1}(t)<b}^{N^u(t)}   \1_{\{X_j^u(t) \geq \frac{3}{2}\ln t +k' \}} e^{-\gamma X_{j}^u(t)+i \sqrt{2} \beta \bar{Y}_{j}^u(t)     }  \right | \geq \epsilon  \right ) =0.
\end{equation*}
\end{lemma}

Notice that we can write:
\begin{align*}
&\sum_{u \in \mathcal{H}_l(t)}  t^{\frac{3 \gamma}{2}} e^{-\gamma X_1^u(t)+i \sqrt{2} \beta \bar{Y}_1^u(t)   } \1_{\{X_1^u(t) \leq \frac{3}{2}\ln t +k \}}  \sum_{j=1, \: t-\tau^u_{j,1}(t)<b}^{N^u(t)}   \1_{\{X_j^u(t) \geq \frac{3}{2}\ln t +k' \}} e^{-\gamma (X_{j}^u(t)-X_{1}^u(t))+i \sqrt{2} \beta (\bar{Y}_{j}^u(t)-\bar{Y}_{1}^u(t)  ) }
\\
&= \sum_{i=1}^{N(t)}  t^{\frac{3 \gamma}{2}} e^{-\gamma X_i(t)+i \sqrt{2} \beta \bar{Y}_i(t)   } \1_{\{ X_i(t)\geq\frac{3}{2}\ln t+k',\,\exists u\in \mathcal{H}_l(t),\, i\in BBM(u),\, X_1^u(t)\leq \frac{3}{2}\ln t+k,\, t-\tau_{i,1}^{u} < b  \}}
\\
&= \sum_{i=1}^{N(t)}  t^{\frac{3 \gamma}{2}} e^{-\gamma X_i(t)+i \sqrt{2} \beta \bar{Y}_i(t)   } \1_{\{ X_i(t)\geq\frac{3}{2}\ln t+k'\}} \1_{\{i \in A(l,t,b,k) \}}
\end{align*}
where $BBM(u)$ is the branching Brownian motion rooted at $u$, $\tau_{i,1^u}$ is the splitting time of $X_i(t)$ and $X_1^u(t)$ and $i\in A(l,t,b,k)$ means:
\begin{equation*}
 \exists u\in \mathcal{H}_l(t),\, i\in BBM(u),\, X_1^u(t)\leq \frac{3}{2}\ln t+k,\, t-\tau_{i,1}^{u}< b.
\end{equation*}
An important observation is that: $\{i\in A(l,t,b,k)\} $ is a set which belong to the sigma field generated by the real branching Brownian motion, and therefore which is indenpendent of the $(\bar{Y}_i(s))_{s\geq 0,\, i \in [1,N(t)]}$. We want to bound 
\begin{equation}
\label{dataffface}\P \left( \left| \sum_{i=1}^{N(t)}  t^{\frac{3 \gamma}{2}} e^{-\gamma X_i(t)+i \sqrt{2} \beta \bar{Y}_i(t)   } \1_{\{ X_i(t)\geq\frac{3}{2}\ln t+k'\}} \1_{\{i \in A(l,t,b,k) \}} \right|\geq \epsilon \right ) .
\end{equation}
Reasoning as in (\ref{stud}), we have:
\begin{align*}
& \P \left( \left| \sum_{i=1}^{N(t)}  t^{\frac{3 \gamma}{2}} e^{-\gamma X_i(t)+i \sqrt{2} \beta \bar{Y}_i(t)   } \1_{\{ X_i(t)\geq\frac{3}{2}\ln t+k'\}} \1_{\{i \in A(l,t,b,k) \}} \right|\geq \epsilon \right ) \\
& \leq e^{-k_1}  + \P\left(\left| \sum_{i=1}^{N(t)}  \1_{\{   \inf_{s\leq t} X_i(s)\geq -k_1,\, \inf_{s\in [t/2,t]}X_i(s)\geq \frac{3}{2} \ln t-L,\, X_i(t)\geq \frac{3}{2} \ln t+k',\, i \in A(l,t,b,k) \}} e^{ -\gamma( X_i(t)-\frac{3}{2}\ln t) +i\beta\sqrt{2}\: \bar{Y}_i(t)}\right| \geq \epsilon \right)  \nonumber
\\
&+ \sum_{L'=L+1}^{2\ln t}\P\left( \exists i \in [|1,N(t)|],\, i \in  \underset{k \in \mathbb{Z}}{\bigcup} \: \underset{v\in \{ t-a_L,...,t\}}{\bigcup}\mathcal{Z}^{k_1,k}(v,L')\right)  \nonumber
\\
&+\P\left(\sum_{L'=L+1}^{2\ln t}\left| \sum_{i=1}^{N(t)}  \1_{\{i \in  \underset{k \in \mathbb{Z}}{\bigcup} \: \underset{v\in \{ t/2,...,t-a_L\}}{\bigcup}\mathcal{Z}^{k_1,k}(v,L'),\, i \in A(l,t,b,k)\}} e^{ -\gamma( X_i(t)-\frac{3}{2}\ln t) +i\beta\sqrt{2}\: \bar{Y}_i(t)}\right| \geq \epsilon \right) 
\\
&\leq e^{-k_1}+ C(1+k_1) e^{-\delta L}+ C e^{k_1}\ln t \: e^{-\delta t} +
\\
&+\frac{1}{\epsilon^2} \E\left( \left|\sum_{i=1}^{N(t)}  \1_{\{   \inf_{s\leq t} X_i(s)\geq -k_1,\, \inf_{s\in [t/2,t]}X_i(s)\geq \frac{3}{2} \ln t-L,\, X_i(t)\geq \frac{3}{2} \ln t+k',\, i \in A(l,t,b,k) \}} e^{ -\gamma( X_i(t)-\frac{3}{2}\ln t) +i\beta\sqrt{2}\: \bar{Y}_i(t)}\right|^2 \right) 
\\
&+  \frac{1}{\epsilon^{{\red \kappa}}} \sum_{L'=L+1}^{2\ln t}\E\left( \left| \sum_{i=1}^{N(t)}  \1_{\{i \in  \underset{k \in \mathbb{Z}}{\bigcup} \: \underset{v\in \{ t/2,...,t-a_L\}}{\bigcup}\mathcal{Z}^{k_1,k}(v,L'),\, i \in A(l,t,b,k)\}} e^{ -\gamma( X_i(t)-\frac{3}{2}\ln t) +i\beta\sqrt{2}\: \bar{Y}_i(t)}\right|^{{\red\kappa}} \right).
\end{align*}
By Jensen's inequality, we therefore get
\begin{align*}
 &e^{-k_1}+ C(1+k_1) e^{-\delta L}+ C e^{k_1}\ln t \:  e^{-\delta t} +
\\
&+ \epsilon^{-2}\E (  \sum_{l=1}^{t}   e^{-2 \beta^2 (t-l)}    \sum_{\tau \in [l,l+1]} e^{ -2\gamma( X_{n_\tau}(\tau)-\frac{3}{2}\ln t)} \sum_{i,j; \: \tau_{i,j}=\tau }  
\\
&  \times \1_{\{   \inf_{s\leq t} X_i(s)\geq -k_1,\, \inf_{s\in [t/2,t]}X_i(s)\geq \frac{3}{2}\ln t-L,\, X_i(t)\geq \frac{3}{2}\ln t+k',\,  i \in A(l,t,b,k) \}} e^{ -\gamma( X_i(t)+X_j(t)-2 X_{n_\tau}(\tau))}  )
\\
&+   \epsilon^{-{\kappa}}\sum_{L=L_0+1}^{2\ln t} \E\left(  \left [ \sum_{l=1}^{t}   e^{-2 \beta^2 (t-l)}    \sum_{v=t/2}^{t-a_L} \sum_{\tau \in [l,l+1]} e^{ -2\gamma( X_{n_\tau}(\tau)-\frac{3}{2}\ln t)} \sum_{i,j; \: \tau_{i,j}=\tau } \1_{\{i \in  \underset{k \in \mathbb{Z}}{\bigcup}\mathcal{Z}^{k_1,k}(v,L),\, i \in A(l,t,b,k)\}}\right. \right.
\\
&\left.\left. e^{ -\gamma( X_i(t)+X_j(t)-2 X_{n_\tau}(\tau))} \right  ] ^{\kappa}  \right)^{ \frac{1}{2}}
\\
&\leq e^{-k_1}+ C(1+k_1) e^{-\delta L}+ C e^{k_1}\ln t \: e^{-\delta t}+ A +  \epsilon^{-{\kappa}}\sum_{L=L_0+1}^{2\ln t}(B_L)^\frac{1}{2}, 
\end{align*}
where $A$ and $B_L$ are the expectations defined respectively in (\ref{proproprop}) (with $k'$ in place of $k$) and (\ref{popopop}) (to get the last inequality, it suffices to remove the indicator $\1_{  i \in A(l,t,b,k)   }$). We then conclude along the same lines as the proofs of lemmas \ref{intermediaire} and \ref{complextension}.

\qed

\subsubsection*{Study of the limit of the BBM} \label{Studylimit}

Here, we define the limit $Z^{(u)}$ defined formally by 
\begin{equation*}
Z^{(u)}= 1+\sum_{j \geq 1} e^{-\gamma \Gamma^{(u)}(\tau_j^{(u)}) -i \sqrt{2} \beta B^{(u)}(\tau_j^{(u)})}   <  e^{-\gamma X+i  \sqrt{2} \beta Y},  \bar{\mathcal{N}}^{(u,j)}_{\Gamma^{(u)}(\tau_j^{(u)}),\tau_j^{(u)}}(dX,dY) >.
\end{equation*} 
 We introduce for all $m \geq 1$
\begin{equation*}
Z^{(u)}_m= 1+\sum_{j=1}^m e^{-\gamma \Gamma^{(u)}(\tau_j^{(u)})- i \sqrt{2} \beta B^{(u)}(\tau_j^{(u)})}   < e^{-\gamma X+i  \sqrt{2} \beta Y},  \bar{\mathcal{N}}^{(u,j)}_{\Gamma^{(u)}(\tau_j^{(u)}),\tau_j^{(u)}}(dX,dY) >
\end{equation*} 
%

We have the following convergence theorem:

\begin{lemma}\label{exislimit}
Let $\beta >0$ and $\gamma>\frac{1}{2}$ be such that $\beta >(1-\gamma)_+$. The sequence $(Z^{(u)}_m)_m$ converges absolutely and almost surely towards a non trivial random variable that we denote
\begin{equation*}
Z^{(u)}= 1+\sum_{j \geq 1} e^{-\gamma \Gamma^{(u)}(\tau_j^{(u)})- i \sqrt{2} \beta B^{(u)}(\tau_j^{(u)})}   < e^{-\gamma X+i  \sqrt{2} \beta Y},  \bar{\mathcal{N}}^{(u,j)}_{\Gamma^{(u)}(\tau_j^{(u)}),\tau_j^{(u)}}(dX,dY) >
\end{equation*}
\end{lemma}

\proof
We consider the less obvious case, i.e. $\beta >0$ and $\gamma \in ]\frac{1}{2},1]$ such that $\gamma +\beta >1$.  
Since the law of $Z^{(u)}$ does not depend on $u$, we consider the case $u=1$ and remove the superscript $^{(u)}$ for clarity. We denote $\E[ . | \Gamma, \tau]$ the conditional expectation with respect to $\Gamma$ and $(\tau_j)_{j \geq 1}$. We introduce an i.i.d. sequence $(X^j)_{j \geq 1}$ of BBMs of law given by the $(X_i(t))_{1 \leq i \leq N(t)}$ of the section Setup and main result.

Now, we have (below $C(\Gamma,(\tau_j)_{j \geq 1})$ denotes a finite constant depending on $\Gamma$ and $(\tau_j)_{j \geq 1}$) 
\begin{align*}
& \E[  \sum_{j=1}^\infty e^{-\gamma \Gamma (\tau_j)}  | < e^{-\gamma X+i  \sqrt{2} \beta Y},  \bar{\mathcal{N}}^{j}_{\Gamma(\tau_j),\tau_j}(dX,dY) > |   | \Gamma, \tau ]  \\
& =  \sum_{j=1}^\infty e^{-\gamma \Gamma(\tau_j)} \E  [  | < e^{-\gamma X+i  \sqrt{2} \beta Y},  \bar{\mathcal{N}}^{j}_{\Gamma(\tau_j),\tau_j}(dX,dY) > |   | \Gamma, \tau ]\\
&  \leq \sum_{j=1}^\infty  \frac{e^{-\gamma \Gamma(\tau_j)}   }{\P( X_1^j(\tau_j) >- \Gamma(\tau_j) )  } \E[ ( \sum_{k,k'}  e^{- \gamma X_k^j(\tau_j) - \gamma X_{k'}^j(\tau_{j}) -2 \beta^2 \tau_{k,k'}^{j,j}}    )^{1/ 2}  |  \Gamma, \tau  ] \\
&  \leq C(\Gamma,(\tau_j)_{j \geq 1}) \sum_{j=1}^\infty e^{-\gamma \Gamma(\tau_j)}  \E[ ( \sum_{k,k'}  e^{- \gamma X_k^j(\tau_j) - \gamma X_{k'}^j(\tau_{j}) -2 \beta^2 \tau_{k,k'}^{j,j}}    )^{1/ 2}  |  \Gamma, \tau  ] \\
&  \leq C(\Gamma,(\tau_j)_{j \geq 1}) \sum_{j=1}^\infty e^{-\gamma \Gamma(\tau_j)}  \E[ ( \sum_{k,k'}  e^{- \gamma X_k^j(\tau_j) - \gamma X_{k'}^j(\tau_{j}) -2 \beta^2 \tau_{k,k'}^{j,j}}    )^{1/ (2\gamma)}  |  \Gamma, \tau  ]^{\gamma} \\
&  \leq C(\Gamma,(\tau_j)_{j \geq 1}) \sum_{j=1}^\infty e^{-\gamma \Gamma(\tau_j)}  \\
& < \infty  
\end{align*}
since $\gamma \in ]1/2,1]$ and where we have used lemma \ref{lemmamoment}.

Finally, it is not hard to see that the limit variable $Z=\underset{m \to \infty}{\lim}Z_m$ is non trivial. Observe that, conditionally to $\Gamma$ and $(\tau_j)_{j \geq 1}$, the variables $(Z_m-Z_{m-1})_{m \geq 1}$ are independent and non constant hence $Z$ is non trivial.

\qed

\begin{lemma}\label{exislimitbis}
Let $\beta >0$ and $\gamma>\frac{1}{2}$ be such that $\beta >(1-\gamma)_+$. The sequence $Z^{(\theta,u)}$ defined by
\begin{equation*}
Z^{(\theta,u)}= 1+\sum_{j \geq 1} e^{-\gamma \Gamma^{(u)}(\tau_j^{(u)})- i \sqrt{2} \beta B^{(u)}(\tau_j^{(u)})}   < \1_{\{X \leq \theta\} }e^{-\gamma X+i  \sqrt{2} \beta Y},  \bar{\mathcal{N}}^{(u,j)}_{\Gamma^{(u)}(\tau_j^{(u)}),\tau_j^{(u)}}(dX,dY) >
\end{equation*}
converges in probability to $Z^{(u)}$ as $\theta$ goes to infinity.
\end{lemma}

\proof
We consider the less obvious case, i.e. $\beta >0$ and $\gamma \in ]\frac{1}{2},1]$ such that $\gamma +\beta >1$. 
Once again, since the law of $Z^{(u)}$ does not depend on $u$, we consider the case $u=1$ and remove the superscript $^{(u)}$ for clarity. We introduce an i.i.d. sequence $(X^j)_{j \geq 1}$ of BBMs of law given by the $(X_i(t))_{1 \leq i \leq N(t)}$ of the section Setup and main result.  Now, we have (below $C(\Gamma,(\tau_j)_{j \geq 1})$ denotes a finite constant depending on $\Gamma$ and $(\tau_j)_{j \geq 1}$)
\begin{align*}
& \E[ |   Z^{(\theta)}-Z  | \:   | \Gamma, \tau ]  \\
& \leq  \E[  \sum_{j=1}^\infty e^{-\gamma \Gamma (\tau_j)}  | < \1_{\{X >\theta\}}e^{-\gamma X+i  \sqrt{2} \beta Y},  \bar{\mathcal{N}}^{j}_{\Gamma(\tau_j),\tau_j}(dX,dY) > |   | \Gamma, \tau ]  \\
& =  \sum_{j=1}^\infty e^{-\gamma \Gamma(\tau_j)} \E  [  | <  \1_{\{X >\theta\}} e^{-\gamma X+i  \sqrt{2} \beta Y},  \bar{\mathcal{N}}^{j}_{\Gamma(\tau_j),\tau_j}(dX,dY) > |   | \Gamma, \tau ]\\
&  \leq C(\Gamma,(\tau_j)_{j \geq 1}) \sum_{j=1}^\infty e^{-\gamma \Gamma(\tau_j)}  \E[ ( \sum_{k,k'}  \1_{\{X_k^j(\tau_j) >\theta\}}  \1_{\{X_{k'}^j(\tau_{j})  >\theta\}} e^{- \gamma X_k^j(\tau_j) - \gamma X_{k'}^j(\tau_{j}) -2 \beta^2 \tau_{k,k'}^{j,j}}    )^{1/ 2}  |  \Gamma, \tau  ].
\end{align*}
Each term $  \E[ ( \sum_{k,k'}  \1_{\{X_k^j(\tau_j) >\theta\}}  \1_{\{X_{k'}^j(\tau_{j})  >\theta\}} e^{- \gamma X_k^j(\tau_j) - \gamma X_{k'}^j(\tau_{j}) -2 \beta^2 \tau_{k,k'}^{j,j}}    )^{1/ 2}  |  \Gamma, \tau  ] $ is bounded by the same quantity without the indicator function and converges to $0$ as $\theta$ goes to infinity. Hence, by the dominated convergence theorem, we have
\begin{equation*}
 \sum_{j=1}^\infty e^{-\gamma \Gamma(\tau_j)}  \E[ ( \sum_{k,k'}  \1_{\{X_k^j(\tau_j) >\theta\}}  \1_{\{X_{k'}^j(\tau_{j})  >\theta\}} e^{- \gamma X_k^j(\tau_j) - \gamma X_{k'}^j(\tau_{j}) -2 \beta^2 \tau_{k,k'}^{j,j}}    )^{1/ 2}  |  \Gamma, \tau  ]  \underset{\theta \to \infty}{\rightarrow} 0.
\end{equation*}
Therefore, the variable $\E[ |   Z^{(\theta)}-Z  | \:   | \Gamma, \tau ]$ converges almost surely to $0$ as $\theta$ goes to infinity. Now, one concludes by using the following inequality for all $\varepsilon>0$
\begin{equation*}
\P(|Z^{(\theta)}-Z  | \geq \varepsilon) \leq \E \left[   \frac{ \E[ |   Z^{(\theta)}-Z  | \:   | \Gamma, \tau ] }{\varepsilon}   \wedge 1 \right ].
\end{equation*}
\qed

\bibliographystyle{plain}

\end{document}